\documentclass[12pt]{article}
\usepackage{e-jc}

\usepackage[german,english]{babel}
\usepackage{amsmath}
\usepackage{amssymb,latexsym}
\usepackage{moreverb}
\usepackage{url}
\usepackage{graphicx}
\usepackage{booktabs}
\usepackage{psfrag}
\usepackage{paralist}
\usepackage{pst-all}

\newtheorem{deff}{Definition}
\newtheorem{lem}[deff]{Lemma}
\newtheorem{prop}[deff]{Proposition}
\newtheorem{thm}[deff]{Theorem}
\newtheorem{cor}[deff]{Corollary}
\newtheorem{conj}[deff]{Conjecture}
\newtheorem{ques}[deff]{Question}
\newtheorem{prob}[deff]{Problem}
\newtheorem{rem}[deff]{Remark}

\makeatletter

\title{$f$-Vectors of $3$-Manifolds}

\author{Frank H.~Lutz\thanks{Supported by the DFG Research Group ``Polyhedral Surfaces'', Berlin}\\
\small Institut f\"ur Mathematik\\[-0.8ex]
\small Technische Universit\"at Berlin\\[-0.8ex]
\small Stra{\ss}e des 17.\ Juni 136, 10623 Berlin, Germany\\
\small \texttt{lutz@math.tu-berlin.de}\\
\and
Thom Sulanke\\
\small Department of Physics\\[-0.8ex]
\small Indiana University\\[-0.8ex]
\small Bloomington, Indiana 47405, USA\\
\small \texttt{tsulanke@indiana.edu}
\and
Ed Swartz\thanks{Paritally supported by NSF grant DMS-0600502}\\
\small Department of Mathematics\\[-0.8ex]
\small Cornell University\\[-0.8ex]
\small Cornell University, Ithaca, NY 14853, USA\\
\small \texttt{ebs@math.cornell.edu}
}

\date{\dateline{May 8, 2008}{May 12, 2009}
      \emph{Dedicated to Anders Bj\"orner on the occasion of his 60th birthday}\\[-15mm]
      \mbox{}}

\begin{document}

\selectlanguage{english}
\maketitle

\bigskip

\begin{abstract}
In 1970, Walkup \cite{Walkup1970} completely described the set of $f$-vectors for the four 
$3$-mani\-folds $S^3$, $S^2\hbox{$\times\hspace{-1.62ex}\_\hspace{-.4ex}\_\hspace{.7ex}$}S^1$,
$S^2\!\times\!S^1$, and ${\mathbb R}{\bf P}^{\,3}$. We improve one of Walkup's main 
restricting inequalities on the set of $f$-vectors of $3$-manifolds. As a consequence 
of a bound by Novik and Swartz~\cite{NovikSwartz2007pre}, we also derive a new lower bound 
on the number of vertices that are needed for a combinatorial $d$-manifold
in terms of its $\beta_1$-coefficient, which partially settles a conjecture of K\"uhnel. 
Enumerative results and a search for small triangulations
with bistellar flips allow us, in combination with the new bounds, to completely determine
the set of $f$-vectors for twenty further $3$-manifolds, that is, for the connected sums 
of sphere bundles $(S^2\!\times\!S^1)^{\# k}$ and twisted sphere bundles 
$(S^2\hbox{$\times\hspace{-1.62ex}\_\hspace{-.4ex}\_\hspace{.7ex}$}S^1)^{\# k}$, 
where $k=2,3,4,5,6,7,8,10,11,14$. For many more $3$-manifolds of different geometric types 
we provide small triangulations and a partial description of their set of $f$-vectors.
Moreover, we show that the $3$-manifold ${\mathbb R}{\bf P}^{\,3}\#\,{\mathbb R}{\bf P}^{\,3}$ 
has (at least) two different minimal $g$-vectors. 
\vspace{15mm}
\mbox{}
\end{abstract}

\section{Introduction}

Let $M$ be a (compact) $3$-manifold (without boundary). According to Moise \cite{Moise1952}, 
$M$ can be triangulated as a (finite) simplicial complex.
If a triangulation of $M$ has \emph{face vector} $f=(f_0,f_1,f_2,f_3)$,
then by Euler's equation, $f_0 - f_1 + f_2 - f_3 = 0$. By double counting
the edges of the triangle-facet incidence graph, $2 f_2 = 4 f_3$.  So
it follows that
\begin{equation}
f=(f_0,f_1,2f_1-2f_0,f_1-f_0).
\end{equation}
In particular,  the number of vertices $f_0$ and the number of edges $f_1$ 
determine the complete $f$-vector of the triangulation.

\begin{thm} {\rm (Walkup \cite{Walkup1970})}\label{thm:Walkup}
For every $3$-manifold $M$ there is a largest integer $\Gamma (M)$
such that
\begin{equation}
f_1\geq 4f_0 - 10 +\Gamma(M)
\end{equation}
for every triangulation of $M$ with $f_0$ vertices and $f_1$ edges 
(with the inequality being tight for at least one triangulation of $M$).
Moreover there is a smallest integer $\Gamma^*(M)\geq \Gamma (M)$
such that for every pair $(f_0,f_1)$ with $f_0\geq 0$ and
\begin{equation}
\binom{f_0}{2}\geq f_1\geq 4f_0 - 10 +\Gamma^*(M)
\end{equation}
there is a triangulation of $M$ with $f_0$ vertices and $f_1$ edges.
Specifically,

\vspace{-2mm}

\begin{itemize}
\addtolength{\itemsep}{-2.5mm}
\item[{\rm (a)}] $\Gamma^*=\Gamma =0$\, for\, $S^3$,
\item[{\rm (b)}] $\Gamma^*=\Gamma =10$\, for\, $S^2\hbox{$\times\hspace{-1.62ex}\_\hspace{-.4ex}\_\hspace{.7ex}$}S^1$,
\item[{\rm (c)}] $\Gamma^*=11$\, and\, $\Gamma =10$\, for\, $S^2\!\times\!S^1$, where, 
                 with the exception $(9,36)$, all pairs $(f_0,f_1)$ with $f_0\geq 0$ 
                 and\, $4f_0\leq f_1\leq\binom{f_0}{2}$ occur,
\item[{\rm (d)}] $\Gamma^*=\Gamma =17$\, for\, ${\mathbb R}{\bf P}^{\,3}$, and
\item[{\rm (e)}] $\Gamma^*(M)\geq\Gamma(M)\geq 18$\, for all other $3$-manifolds $M$.
\addtolength{\itemsep}{2.5mm}
\end{itemize}

\end{thm}
By definition, $\Gamma(M)$ and $\Gamma^*(M)$ are topological invariants of $M$,
with $\Gamma(M)$ determining the range of pairs $(f_0,f_1)$ for which triangulations
of $M$ can occur, whereas $\Gamma^*(M)$ ensures that for all pairs $(f_0,f_1)$
in the respective range there indeed are triangulations with the corresponding $f$-vectors.

\begin{rem}
Walkup originally stated Theorem~\ref{thm:Walkup} in terms of the constants\, $\gamma=-10+\Gamma$ 
and\, $\gamma^*=-10+\Gamma^*$. As we will see in Section~\ref{sec:3d}, our choice 
of the constant\, $\Gamma(M)$ (as well as of\, $\Gamma^*(M)$) is more naturally related to 
the $g_2$-entries of the $g$-vectors of triangulations of a $3$-manifold $M$: 
$\Gamma(M)$ is the smallest $g_2$ that is possible for all triangulations of $M$.
\end{rem}

In the following section, we review some of the basic facts on $f$- and $g$-vectors
of triangulated $d$-manifolds and how they change under (local) modifications 
of the triangulation. Moreover, we derive a new bound on the minimal number
of vertices for a triangulable $d$-manifold depending on its $\beta_1$-coefficient.
In Section~\ref{sec:3d}, we discuss the $f$- and $g$-vectors of $3$-manifolds
in more detail and introduce tight-neighborly triangulations. 
Section~\ref{sec:g2irreducible} is devoted to the proof of an improvement of a bound by Walkup and to the notion of 
$g_2$-irreducible triangulations. In Section~\ref{sec:enumeration}
we completely enumerate all $g_2$-irreducible triangulations of $3$-manifolds 
with $g_2\leq 20$ and all potential $g_2$-irreducible triangulations of $3$-mani\-folds with $f_0 \leq 15$.
Section~\ref{sec:examples} presents small triangulations of different geometric types,
in particular, examples of Seifert manifolds from the six Seifert geometries
as well as triangulations of hyperbolic $3$-manifolds. With the help of these
small triangulations we establish upper bounds on the invariants $\Gamma$
and $\Gamma^*$ of the respective manifolds. For the $3$-manifold 
${\mathbb R}{\bf P}^{\,3}\#\,{\mathbb R}{\bf P}^{\,3}$ we show that it has (at least) 
two different minimal $g$-vectors. Finally, we extend Walkup's Theorem~\ref{thm:Walkup} by
completely characterizing the set of $f$-vectors of the twenty $3$-manifolds $(S^2\!\times\!S^1)^{\# k}$ 
and $(S^2\hbox{$\times\hspace{-1.62ex}\_\hspace{-.4ex}\_\hspace{.7ex}$}S^1)^{\# k}$ with
$k=2,3,4,5,6,7,8,10,11,14$. (In dimensions $d\geq 4$, a complete description 
of the set of $f$-vectors is only known for the six $4$-manifolds $S^4$~\cite{McMullen1971},
$S^3\!\times\!S^1$, ${\mathbb C}{\bf P}^{\,2}$, $K3$-surface, $(S^2\!\times\!S^1)^{\# 2}$ \cite{Swartz2007apre},
and $S^3\hbox{$\times\hspace{-1.62ex}\_\hspace{-.4ex}\_\hspace{.7ex}$}S^1$ \cite{ChestnutSapirSwartz2008}.)
In~Section~\ref{sec:gamma_complexity} we compare the invariant $\Gamma(M)$ to Matveev's complexity measure $c(M)$.

\section{Face Numbers and (Local) Modifications}

Let $K$ be a triangulation of a $d$-manifold $M$ with $f$-vector $f(K)=(f_0(K),\dots,f_d(K))$
(and with $f_{-1}(K)=1$), that is, $f_i(K)$ denotes the number of $i$-dimensional faces of $K$.
For simplicity, we write $f=(f_0,\dots,f_d)$, and we define numbers $h_i$ by 

\vspace{-3.5mm}

\begin{equation}
h_k=\sum^k_{i=0}(-1)^{k-i}\binom{d+1-i\,}{d+1-k}f_{i-1}.\label{eq:hk}
\end{equation}

\vspace{-.5mm}

\noindent
The vector $h=(h_0, \dots, h_{d+1})$ is called the \emph{$h$-vector} of $K$.
Moreover, the $g$-\emph{vector} \mbox{$g=(g_0,\dots ,g_{\lfloor {(d+1)/2} \rfloor})$} of $K$
is defined by $g_0=1$ and $g_k=h_k-h_{k-1}$, for\, $k\ge 1$, which gives

\vspace{-7.5mm}

\begin{eqnarray}
g_k & = & \sum^k_{i=0}(-1)^{k-i}\binom{d+2-i\,}{d+2-k}f_{i-1}.\label{eq:gk}
\end{eqnarray}

\vspace{-.5mm}

\noindent
In particular,

\vspace{-7mm}

\begin{eqnarray}
g_1 & = & f_0-(d+2),\label{eq:g1} \\[-1.5mm]
g_2 & = & f_1-(d+1)f_0+\binom{d+2}{2}. \label{eq:g2}
\end{eqnarray}

\vspace{-.5mm}

Let ${\cal H}^d$ be the class of triangulated $d$-manifolds 
that can be obtained from the boundary of the $(d+1)$-simplex
by a sequence of the following three operations:

\vspace{.5mm}

\begin{itemize}
\addtolength{\itemsep}{-1mm}
\item[$S$] Subdivide a facet with one new vertex in the interior of the facet.
\item[$H$] Form a handle (oriented or non-oriented) by identifying a pair of facets in $K\in{\cal H}^d$
           and removing the interior of the identified facet in such a way that the resulting
           complex is still a simplicial complex (i.e., the distance in the $1$-skeleton of $K$ 
           between every pair of identified vertices must be at least three).
\item[\#]  Form the connected sum of $K_1,K_2\in{\cal H}^d$ by identifying a pair
           of facets, one from each complex, and then removing the interior
           of the identified facet.
\addtolength{\itemsep}{1mm}
\end{itemize}

\pagebreak

\noindent
For the operations $S$, $H$, and \# the resulting triangulations 
depend on the particular choices of the facets and, in the case 
of $H$ and \#, on the respective identifications. 
However, all triangulated $d$-manifolds 
in the class ${\cal H}^d$ are of the following topological types:
the $d$-sphere $S^d$, connected sums $(S^{d-1}\!\times\!S^1)^{\# k}$ of the orientable sphere product
$S^{d-1}\!\times\!S^1$ for $k\geq 1$,
or connected sums $(S^{d-1}\hbox{$\times\hspace{-1.62ex}\_\hspace{-.4ex}\_\hspace{.7ex}$}S^1)^{\# l}$
of the twisted sphere product $S^{d-1}\hbox{$\times\hspace{-1.62ex}\_\hspace{-.4ex}\_\hspace{.7ex}$}S^1$
for $l\geq 1$.

Let $K$, $K_1$, and $K_2$ be arbitrary triangulated $d$-manifolds with $f$-vectors
\begin{eqnarray*}
f(K)   & = & (f_0(K),\dots,f_d(K)),\\
f(K_1) & = & (f_0(K_1),\dots,f_d(K_1)),\\
f(K_2) & = & (f_0(K_2),\dots,f_d(K_2)),
\end{eqnarray*}
and $f_{-1}(K)=f_{-1}(K_1)=f_{-1}(K_2)=1$. Again, let $SK$ be the triangulated $d$-manifold 
obtained from $K$ by performing the subdivision operation $S$ on some facet of $K$, $HK$ 
be the triangulated $d$-manifold obtained from $K$ by performing the handle addition operation~$H$ 
on some (admissible) pair of facets of $K$, and $K_1\#\pm K_2$ be the
triangulated $d$-manifold obtained from $K_1$ and $K_2$ 
by the connected sum operation $\#$ on some pair of facets of $K_1$ and $K_2$.
Then the $f$-vectors of $SK$, $HK$, and $K_1\#\pm K_2$ have entries 
\begin{eqnarray}
f_k(SK)           & = & f_k(K)+\binom{d+1}{k},\quad \mbox{for}\quad 0\leq k\leq d-1,\\[-.5mm]
f_d(SK)           & = & f_d(K)+d,\\[1.5mm]
f_k(HK)           & = & f_k(K)-\binom{d+1}{k+1},\quad\mbox{for}\quad 0\leq k\leq d-1,\\[-.5mm]
f_d(HK)           & = & f_d(K)-2,\\[1.5mm]
f_k(K_1\#\pm K_2) & = & f_k(K_1)+f_k(K_2)-\binom{d+1}{k+1},\quad\mbox{for}\quad 0\leq k\leq d-1,\\[-.5mm]
f_d(K_1\#\pm K_2) & = & f_d(K_1)+f_d(K_2)-2.
\end{eqnarray}
In particular, it follows that
\begin{eqnarray}
  g_1(SK)           & = & g_1(K)+1,\label{eq:g1sm}\\[-.5mm]
  g_k(SK)           & = & g_k(K),\quad\mbox{for}\quad 2\leq k\le \lfloor {(d+1)/2}\rfloor,\\[1.5mm]
  g_1(HK)           & = & g_1(K)-(d+1),\\[-.5mm]
  g_k(HK)           & = & g_k(K)+(-1)^k\binom{d+2}{k},\quad\mbox{for}\quad 2\leq k\le \lfloor {(d+1)/2}\rfloor,\\[1.5mm]
  g_1(K_1\#\pm K_2) & = & g_1(K_1)+g_1(K_2)+1,\\[-.5mm]
  g_k(K_1\#\pm K_2) & = & g_k(K_1)+g_k(K_2),\quad\mbox{for}\quad 2\leq k\le \lfloor {(d+1)/2}\rfloor.\label{eq:gksum}
\end{eqnarray}

\begin{conj} {\rm (Kalai~\cite{Kalai1987})}\label{conj:kalai}
Let $K$ be a connected triangulated $d$-manifold with $d\geq 3$.
Then
\begin{equation}
g_2(K)\geq\binom{d+2}{2}\beta_1(K;{\mathbb Q}).\label{eq:kalai}
\end{equation}
\end{conj}

In \cite{Swartz2007apre}, Swartz verified Kalai's conjecture for all $d\geq 3$ when $\beta_1(K;{\mathbb Q})=1,$ 
and for orientable $K$ when $d \ge 4$ and $\beta_2(K,{\mathbb Q})=0.$

\begin{thm} {\rm (Novik and Swartz~\cite{NovikSwartz2007pre})}\label{thm:NovikSwartz}
Let ${\mathbb K}$ be any field and let $K$ be a (connected) triangulation of a 
${\mathbb K}$-orientable ${\mathbb K}$-homology $d$-dimensional
manifold with $d\geq 3$. Then
\begin{equation}
g_2(K)\geq\binom{d+2}{2}\beta_1(K;{\mathbb K}).\label{eq:NovikSwartz}
\end{equation}
Furthermore, if $g_2=\binom{d+2}{2}\beta_1(K;{\mathbb K})$ and $d\geq 4,$ then $K\in{\cal H}^d$.
\end{thm}
Since any $d$-manifold (without boundary) is orientable over ${\mathbb K}$ if ${\mathbb K}$ 
has characteristic two, and in this case $\beta_1({\mathbb K}) \ge \beta_1({\mathbb Q})$ 
(universal coefficient theorem), this theorem proves Conjecture \ref{conj:kalai}. 

Combining (\ref{eq:NovikSwartz}) and (\ref{eq:g2}) with the obvious
inequality $f_1\leq\binom{f_0}{2}$ yields

\begin{center}
\begin{tabular}{lllll}
$\displaystyle\binom{d+2}{2}\beta_1$ & $\leq$ & $g_2$  & $=$    & $f_1-(d+1)f_0+\displaystyle\binom{d+2}{2}$ \\[5mm]
                                     &        &        & $\leq$ & $\displaystyle\binom{f_0}{2}-(d+1)f_0+\displaystyle\binom{d+2}{2}$
\end{tabular}
\end{center}
or, equivalently,
\begin{eqnarray}
f_0^2-(2d+3)f_0+(d+1)(d+2)(1-\beta_1) & \geq & 0.\label{eq:f0bound}
\end{eqnarray}

\begin{thm}\label{thm:f0m}
Let ${\mathbb K}$ be any field and let $K$ be a ${\mathbb K}$-orientable
triangulated $d$-mani\-fold with $d\geq 3$. Then
\begin{equation}
f_0(K)\geq\Bigl\lceil\tfrac{1}{2}\Bigl((2d+3)+\sqrt{1+4(d+1)(d+2)\beta_1(K;{\mathbb K})}\Bigl)\Bigl\rceil.\label{eq:f0m}
\end{equation}
\end{thm}

Inequality~(\ref{eq:f0bound}) can also be written in the form $\binom{f_0-d-1}{2}\geq \binom{d+2}{2}\beta_1$.
Its proof settles K\"uhnel's conjectured bounds $\binom{f_0-d+j-2}{j+1}\geq \binom{d+2}{j+1}\beta_j$ 
(cf.\ \cite{Lutz2005bpre}, with $1\leq j\leq\lfloor\frac{d-1}{2}\rfloor$) in the cases with $j=1$.

According to Brehm and K\"uhnel \cite{BrehmKuehnel1987},
we further have for all $(j-1)$-connected but not $j$-connected combinatorial $d$-manifolds $K$, with $1\leq j< d/2$, that
\begin{equation}
f_0(K)\geq 2d+4-j\label{eq:f0bk}.
\end{equation}

While the bound (\ref{eq:f0m}) becomes trivial for manifolds with $\beta_1=0$, with the $d$-sphere~$S^d$
admitting triangulations in the full range $f_0(S^d)\geq d+2$, the inequality (\ref{eq:f0bk})
yields stronger restrictions for higher-connected manifolds.
In contrast, for all non-simply connected combinatorial $d$-manifolds $K$ the bound (\ref{eq:f0bk}) uniformly gives
\begin{equation}
f_0(K)\geq 2d+3\label{eq:f0bkj=1},
\end{equation}
whereas the bound (\ref{eq:f0m}) explicitly depends on the $\beta_1$-coefficient.

\pagebreak

In the case $\beta_1=1$, the bounds (\ref{eq:f0m}) and (\ref{eq:f0bk}) coincide with (\ref{eq:f0bkj=1})
and are sharp for 
\begin{itemize}
\addtolength{\itemsep}{-1.5mm}
\item $S^{d-1}\!\times\!S^1$\, if $d$ is even \cite{Kuehnel1986a-series,KuehnelLassmann1996-bundle},
\item $S^{d-1}\hbox{$\times\hspace{-1.62ex}\_\hspace{-.4ex}\_\hspace{.7ex}$}S^1$\, if $d$ is odd \cite{Kuehnel1986a-series,KuehnelLassmann1996-bundle}, 
\addtolength{\itemsep}{1.5mm}
\end{itemize}~while\,  $f_0(S^{d-1}\!\times\!S^1)\geq 2d+4$\, for $d$ odd
and\, $f_0(S^{d-1}\hbox{$\times\hspace{-1.62ex}\_\hspace{-.4ex}\_\hspace{.7ex}$}S^1)\geq 2d+4$\, for $d$ even;
see \cite{BagchiDatta2008,ChestnutSapirSwartz2008}.

If $K$ is a triangulated $2$-manifold with Euler characteristic $\chi(K)$,
then by Heawood's inequality \cite{Heawood1890},
\begin{equation}
f_0\geq\Bigl\lceil\tfrac{1}{2}\Bigl(7+\sqrt{49-24\chi (K)}\Bigl)\Bigl\rceil.\label{eq:heawood}
\end{equation}
For an orientable surface $K$ of genus $g$ the Euler characteristic of $K$ is $2-2g$.
Therefore 
$f_0\geq\lceil\tfrac{1}{2}(7+\sqrt{1+48g})\rceil$,
whereas $\chi(K)=2-u$ for a non-orientable surface $K$ of genus $u$
and hence $f_0\geq\lceil\tfrac{1}{2}(7+\sqrt{1+24u})\rceil$.
These bounds all coincide with 
\begin{equation}
f_0\geq\Bigl\lceil\tfrac{1}{2}\Bigl(7+\sqrt{1+48\tfrac{\beta_1(K;{\mathbb Z}_2)}{2}}\Bigl)\Bigl\rceil
\end{equation}
or, equivalently, $\binom{f_0-3}{2}\geq\binom{4}{2}\tfrac{\beta_1(K;{\mathbb Z}_2)}{2}$, 
where the factor $\frac{1}{2}$ on the right hand side compensates the doubling of homology
in the middle homology of even dimensional manifolds by Poincar\'e duality; see \cite{Lutz2005bpre}
for K\"uhnel's conjectured higherdimensional analogues of this bound.

Heawood's bound (\ref{eq:heawood}) is sharp, except in the cases of the orientable surface of genus~$2$, 
the Klein bottle, and the non-orientable surface of genus~$3$. Each of these requires
an extra vertex  to be added. The construction of series of examples of 
vertex-minimal triangulations was completed in 1955 
for all non-orientable surfaces by Ringel~\cite{Ringel1955} and 
in 1980 for all orientable surfaces by Jungerman and Ringel~\cite{JungermanRingel1980}.

\begin{ques}\label{ques:connectedsums}
Is inequality~(\ref{eq:f0m}) sharp for all but finitely many connected sums $(S^{d-1}\!\times\!S^1)^{\# k}$ 
of sphere products as well as for all but finitely many connected sums 
$(S^{d-1}\hbox{$\times\hspace{-1.62ex}\_\hspace{-.4ex}\_\hspace{.7ex}$}S^1)^{\# k}$
of twisted sphere products in every fixed dimension $d\geq 3$?
\end{ques}

\begin{prob}\label{prob:connectedsums}
Construct series of vertex-minimal triangulations of $(S^{d-1}\!\times\!S^1)^{\# k}$
and of $(S^{d-1}\hbox{$\times\hspace{-1.62ex}\_\hspace{-.4ex}\_\hspace{.7ex}$}S^1)^{\# k}$ for $d\geq 3$.
Can the examples be chosen to lie in the class ${\cal H}^d$?
\end{prob}
The only known series of such vertex-minimal triangulations \label{ques:connectedsum}
are the ones mentioned above of the $d$-sphere $S^d\cong (S^{d-1}\!\times\!S^1)^{\# 0}\cong (S^{d-1}\hbox{$\times\hspace{-1.62ex}\_\hspace{-.4ex}\_\hspace{.7ex}$}S^1)^{\# 0}$, 
triangulated as the boundary of the $(d+1)$-simplex with $d+2$ vertices, and for $k=1$
the vertex-minimal triangulations of $(S^{d-1}\!\times\!S^1)^{\# 1}$ 
and $(S^{d-1}\hbox{$\times\hspace{-1.62ex}\_\hspace{-.4ex}\_\hspace{.7ex}$}S^1)^{\# 1}$.

A first sporadic vertex-minimal $4$-dimensional example with $k=3$ was recently discovered by Bagchi and Datta~\cite{BagchiDatta2008bpre}.  
They construct a triangulation of $(S^3\hbox{$\times\hspace{-1.62ex}\_\hspace{-.4ex}\_\hspace{.7ex}$}S^1)^{\# 3}$ 
in ${\cal H}^d$ with $15$ vertices and $g_2=45$, both the minimums required by (\ref{eq:NovikSwartz}) and (\ref{eq:f0m}).
For $3$-dimensional examples with $k=2,3,4,5,6,7,8,10,11,14$, see Theorem~\ref{thm:extensionWalkup} below.

\pagebreak

Besides subdivisions, handle additions, and connected sums,
\emph{bistellar flips} (also called Pachner moves \cite{Pachner1986}) 
are a very useful class of local modifications of triangulations.
\begin{deff} {\rm \cite{Pachner1986}}
Let $K$ be a triangulated $d$-manifold. If $A$ is a\, $(d-i)$-face of~$K$, $0\leq i\leq d$,
such that the link of $A$ in $K$, ${\rm Lk}\,A$, 
is the boundary\, $\partial (B)$ of an $i$-simplex $B$ 
that is not a face of $K$, then the operation $\Phi_A$ on $K$
defined by
$$\Phi_A(K):=(K\backslash (A*\partial (B)))\cup (\partial (A)*B)$$
is a\, \emph{bistellar $i$-move} (with $*$ the join operation for simplicial complexes).
\end{deff}
In particular, the subdivision operation $S$ from above on any facet $A$ of $K$
coincides with the bistellar $0$-move on this facet.

If\, $K'$ is obtained from $K$ by a bistellar $i$-move, $0\le i\le \lfloor {(d-1)/2} \rfloor$, 
then
\begin{eqnarray}
g_{i+1}(K') & = & g_{i+1}(K)+1\\
g_{k}(K')   & = & g_{k}(K)\, \quad\mbox{for all}\quad k\neq i+1.
\end{eqnarray}
If\, $d$ is even and\, $i=\frac{d}{2}$, 
then
\begin{eqnarray}
 g_{k}(K') & = & g_{k}(K)\, \quad\mbox{for all}\quad k.
\end{eqnarray}

Bistellar flips  can be used to navigate through
the set of triangulations of a $d$-manifold,
 with the objective of obtaining a small, or perhaps
even vertex-minimal, triangulation of this manifold. 
A simulated annealing type strategy for this aim
is described in \cite{BjoernerLutz2000}. The reference
 \cite{BjoernerLutz2000} also contains further background
on combinatorial topology aspects of bistellar flips.
A basic implementation of the bistellar flip heuristics
is \cite{Lutz_BISTELLAR}. The \texttt{bistellar} client
of the \texttt{polymake}-system~\cite{polymake} allows for fast computations
for rather large triangulations, as we will need in Section~\ref{sec:examples}.

\section{Face Numbers and (Local) Modifications\\ for 3-Manifolds}
\label{sec:3d}

Let $K$ be a triangulated $3$-manifold with $f$-vector $f=(f_0,f_1,2f_1-2f_0,f_1-f_0)$.
The relations~(\ref{eq:g1}), (\ref{eq:g2}), (\ref{eq:g1sm})--(\ref{eq:gksum}) then read,
\begin{eqnarray}
g_1 & = & f_0-5, \\
g_2 & = & f_1-4f_0+10\label{eq:g23d}
\end{eqnarray}
and
\begin{eqnarray}
g_1(SK) & = & g_1(K)+1, \\
g_2(SK) & = & g_2(K),\label{eq:g2sm3d}\\[2mm]
g_1(HK) & = & g_1(K)-4, \\
g_2(HK) & = & g_2(K)+10,\label{eq:g2hm3d}\\[2mm]
g_1(K_1\#\pm K_2) & = & g_1(K_1)+g_1(K_2)+1, \\
g_2(K_1\#\pm K_2) & = & g_2(K_1)+g_2(K_2).\label{eq:g2m1m2}
\end{eqnarray}
For a $3$-manifold $M$, $\Gamma(M)$ is the smallest $g_2$ that is possible 
for all triangulations of $M$. Hence, the following lemma follows immediately from (\ref{eq:g2m1m2}) and Theorem~\ref{thm:Walkup}. 

\begin{lem}\label{lem:handleadd}
Let $M$, $M_1$, and $M_2$ be $3$-manifolds. Then
\begin{eqnarray}
\Gamma(M_1\#\pm M_2) & \leq & \Gamma(M_1)+\Gamma(M_2),\label{eq:gammasum}\\
\Gamma(M\#(S^2\!\times\!S^1)^{\# k}) & \leq & \Gamma(M)+10k,\label{eq:gammaorhandle}\\
\Gamma(M\#(S^2\hbox{$\times\hspace{-1.62ex}\_\hspace{-.4ex}\_\hspace{.7ex}$}S^1)^{\# k}) & \leq & \Gamma(M)+10k.\label{eq:gammanorhandle}
\end{eqnarray}
\end{lem}

As a consequence of Theorem~\ref{thm:NovikSwartz} and Lemma~\ref{lem:handleadd}:

\begin{cor}
For every\, $k\in{\mathbb N}$,
\begin{equation}
\Gamma((S^2\!\times\!S^1)^{\# k}) = \Gamma(S^2\hbox{$\times\hspace{-1.62ex}\_\hspace{-.4ex}\_\hspace{.7ex}$}S^1)^{\# k}) = 10k.
\end{equation}
\end{cor}

\begin{conj}  \label{conj:connectedsum}
Let $M$, $M_1$, and $M_2$ be $3$-manifolds. Then
\begin{eqnarray}
\Gamma(M_1\#\pm M_2) & = & \Gamma(M_1)+\Gamma(M_2),\\
\Gamma(M\#(S^2\!\times\!S^1)^{\# k}) & = & \Gamma(M)+10k,\\
\Gamma(M\#(S^2\hbox{$\times\hspace{-1.62ex}\_\hspace{-.4ex}\_\hspace{.7ex}$}S^1)^{\# k}) & = & \Gamma(M)+10k.
\end{eqnarray}
\end{conj}

\noindent 
While the latter two equalities above would follow from the first, 
it may be the case that only these two special cases hold. 

By Theorem~\ref{thm:f0m}, inequality (\ref{eq:f0m})
holds for all ${\mathbb K}$-orientable triangulated $3$-mani\-folds $K$,
that is,
\begin{equation}
f_0(K)\geq\Bigl\lceil\tfrac{1}{2}\Bigl(9+\sqrt{1+80\beta_1(K;{\mathbb K})}\Bigl)\Bigl\rceil.\label{eq:f0m3d}
\end{equation}

We next consider the class of connected sums
$(S^2\!\times\!S^1)^{\# k}$ and $(S^2\hbox{$\times\hspace{-1.62ex}\_\hspace{-.4ex}\_\hspace{.7ex}$}S^1)^{\# k}$
with $\beta_1=k$ for $k\in{\mathbb N}$. For this class, inequality (\ref{eq:f0bound})
can be interpreted as an upper bound on the number $k$ 
for which the corresponding connected sums can have triangulations with $f_0$ vertices. 
If inequality (\ref{eq:f0bound}) is sharp, then $f_1=\binom{f_0}{2}$, i.e., such a triangulation 
must be \emph{neighborly} with complete $1$-skeleton.
We therefore call triangulations of connected sums of the sphere bundles
$S^2\!\times\!S^1$ and $S^2\hbox{$\times\hspace{-1.62ex}\_\hspace{-.4ex}\_\hspace{.7ex}$}S^1$
for which inequality (\ref{eq:f0bound}) is tight \emph{tight-neighborly}.
In the case of equality,
\begin{equation}
\frac{(f_0-9)f_0}{20}=k-1,\label{eq:f0eq}
\end{equation}
the right hand side of (\ref{eq:f0eq}) is integer, 
and therefore, the left hand side is integer as well.
This is possible if and only if
\begin{equation}
f_0\equiv 0,4,5,9\,\mbox{mod}\,20,
\end{equation}
with the additional requirement that $f_0\geq 5$. 
\begin{table}
\small\centering
\defaultaddspace=0.2em
\caption{Parameters for tight-neighborly triangulations}\label{tbl:tight_neighborly}
\begin{tabular}{@{}l@{\hspace{12mm}}l@{}}
\\\toprule
 \addlinespace
 \addlinespace
    $f_0$         &  $k$     \\ \midrule
 \addlinespace
 \addlinespace
 \addlinespace
    $20m$         &  $20m^2-9m+1$ \\
 \addlinespace
    $4+20m$       &  $20m^2-m$ \\
 \addlinespace
    $5+20m$       &  $20m^2+m$ \\
 \addlinespace
    $9+20m$       &  $20m^2+9m+1$ \\
 \addlinespace
 \addlinespace
\bottomrule
\end{tabular}
\end{table}
Table~\ref{tbl:tight_neighborly} gives the possible parameters $(f_0,k)$ for 
tight-neighborly triangulations. The first two pairs are $(f_0,k)=(5,0)$
and $(f_0,k)=(9,1)$, for which we have the triangulation of $S^3$ as the
boundary $\partial\Delta^4$ of the $4$-simplex $\Delta^4$ and Walkup's unique $9$-vertex triangulation \cite{Walkup1970}
of $S^2\hbox{$\times\hspace{-1.62ex}\_\hspace{-.4ex}\_\hspace{.7ex}$}S^1$, respectively.
There is no triangulation of $S^2\!\times\!S^1$ with $9$ vertices.

\begin{ques}
Are there $3$-dimensional tight-neighborly triangulations for $k>1$?
\end{ques}
The first two cases would be $(f_0,k)=(20,12)$ and $(f_0,k)=(24,19)$.

Tight-neighborly triangulations are possible candidates for ``tight triangulations''
in the following sense (cf., \cite{Kuehnel1995-book,KuehnelLutz1999}):
A simplicial complex $K$ with vertex-set $V$
is \emph{tight} if for any subset $W\subseteq V$ of vertices the induced homomorphism
$$H_\ast (\langle W \rangle \cap K;{\mathbb K}) \to H_\ast (K;{\mathbb K})$$
is injective, where $\langle W \rangle$ denotes the face of the $(|V|-1)$-simplex $\Delta^{|V|-1}$
spanned by $W$.

Obviously, we can extend the concept of tight-neighborly triangulations to any dimension $d\geq 2$:
Triangulations of connected sums of sphere bundles $S^{(d-1)}\!\times\!S^1$ 
and $S^{(d-1)}\hbox{$\times\hspace{-1.62ex}\_\hspace{-.4ex}\_\hspace{.7ex}$}S^1$
are \emph{tight-neighborly} if inequality (\ref{eq:f0bound}) is tight.

By Theorem~\ref{thm:NovikSwartz}, every triangulation $K$ of a ${\mathbb K}$-orientable 
${\mathbb K}$-homology $d$-dimensional manifold with $d\geq 4$ for which (\ref{eq:f0bound}) is tight
lies in ${\cal H}^d$ and therefore is a tight-neighborly connected sum of sphere bundles $S^{(d-1)}\!\times\!S^1$
or $S^{(d-1)}\hbox{$\times\hspace{-1.62ex}\_\hspace{-.4ex}\_\hspace{.7ex}$}S^1$.

\begin{conj}\label{conj:tight_neighborly}
Tight-neighborly triangulations are tight.
\end{conj}
The conjecture holds for surfaces (i.e., for $d=2$) \cite[Sec.~2D]{Kuehnel1995-book},
for $k=0$ (that is, for the triangulation of $S^d$ as the boundary of the $(d+1)$-simplex) \cite[Sec.~3A]{Kuehnel1995-book}, 
and for $k=1$, in which case there is a unique and tight triangulation with $2d+3$ vertices in every dimension $d\geq 2$ 
(see \cite{Kuehnel1986a-series,Moebius1886,Walkup1970} for existence, \cite{BagchiDatta2008,ChestnutSapirSwartz2008} 
for uniqueness, and \cite[Sec.~5B]{Kuehnel1995-book} for tightness). 

For the sporadic Bagchi-Datta example~\cite{BagchiDatta2008bpre} we used the computational methods 
from \cite{KuehnelLutz1999} to determine the tightness.

\begin{prop}
The tight-neighborly $4$-dimensional $15$-vertex example of Bagchi and Datta 
with $k=3$ is tight.
\end{prop}

Most recently, Conjecture~\ref{conj:tight_neighborly} was settled in even dimensions $d\geq 4$ 
by Effenberger~\cite{Effenberger2009pre}. In particular, this also yields the tightness 
of the Bagchi-Datta example.

\pagebreak

\section{\mathversion{bold}$g_2$-Irreducible\mathversion{normal} Triangulations}
\label{sec:g2irreducible}

The main idea behind the proof of Theorem \ref{thm:Walkup} is that triangulations which minimize $g_2$  
have several special combinatorial properties. A {\it missing facet} of a triangulated $d$-manifold $K$ 
is a subset $\sigma$ of the vertex set of cardinality $d+1$ such that $\sigma \notin K,$ 
but every proper subset of $\sigma$ is a face of $K.$  

\begin{deff}
Let $K$ be a triangulation of a $3$-manifold $M$. Then $K$ is \emph{$g_2$-mini\-mal} 
if $g_2(K') \ge g_2(K)$ for all other triangulations $K'$ of $M$,
i.e., $g_2(K)=\Gamma(M)$. The triangulation $K$ is \emph{$g_2$-irreducible}  if
the following hold:
\begin{enumerate}
\addtolength{\itemsep}{-1mm}
\item $K$ is $g_2$-minimal.
\item $K$ is not the boundary of the $4$-simplex.
\item $K$ does not have any missing facets.
\addtolength{\itemsep}{1mm}
\end{enumerate}
\end{deff}

The reason for introducing the third condition is the following folk theorem.  For a complete proof, see \cite[Lemma 1.3]{BagchiDatta2008}.

\begin{thm} \label{thm:missingfacet}
Let $K$ be a triangulated $3$-manifold.  Then $K$ has a missing facet if and only if $K$ equals $K_1 \# K_2$ or $HK^\prime.$
\end{thm}

So, a triangulation $K$ which realizes the minimum $g_2$ for a particular $3$-manifold $M$ is either $g_2$-irreducible, 
or is of the form $K_1 \# K_2$ or $HK^\prime$, where the component triangulations realize their minimum $g_2.$ 
The remainder of this section is devoted to proving the following.

\begin{thm}  \label{thm:newwalkup}
If $K$ is $g_2$-irreducible, then 
\begin{equation}\label{eq:newwalkup}
f_1(K)- \tfrac{9}{2} f_0(K) > \tfrac{1}{2}.
\end{equation}
\end{thm}

Walkup originally proved that for $g_2$-irreducible $K$, $f_1(K)-\tfrac{9}{2} f_0(K) > 0.$  All that is needed to get the slight improvement 
we require is a little more care.  Walkup's original result plus Theorem \ref{thm:missingfacet} are already enough 
to prove that for a fixed $\Gamma$ there are only finitely many $3$-manifolds such that $\Gamma(M) \le \Gamma$ \cite{Swartz2008}
(see also the next section).  
With the exception of Theorem \ref{thm:newwalkup}, all of the remaining results in this section first appeared in \cite{Walkup1970} 
and we refer the reader there for the proofs.

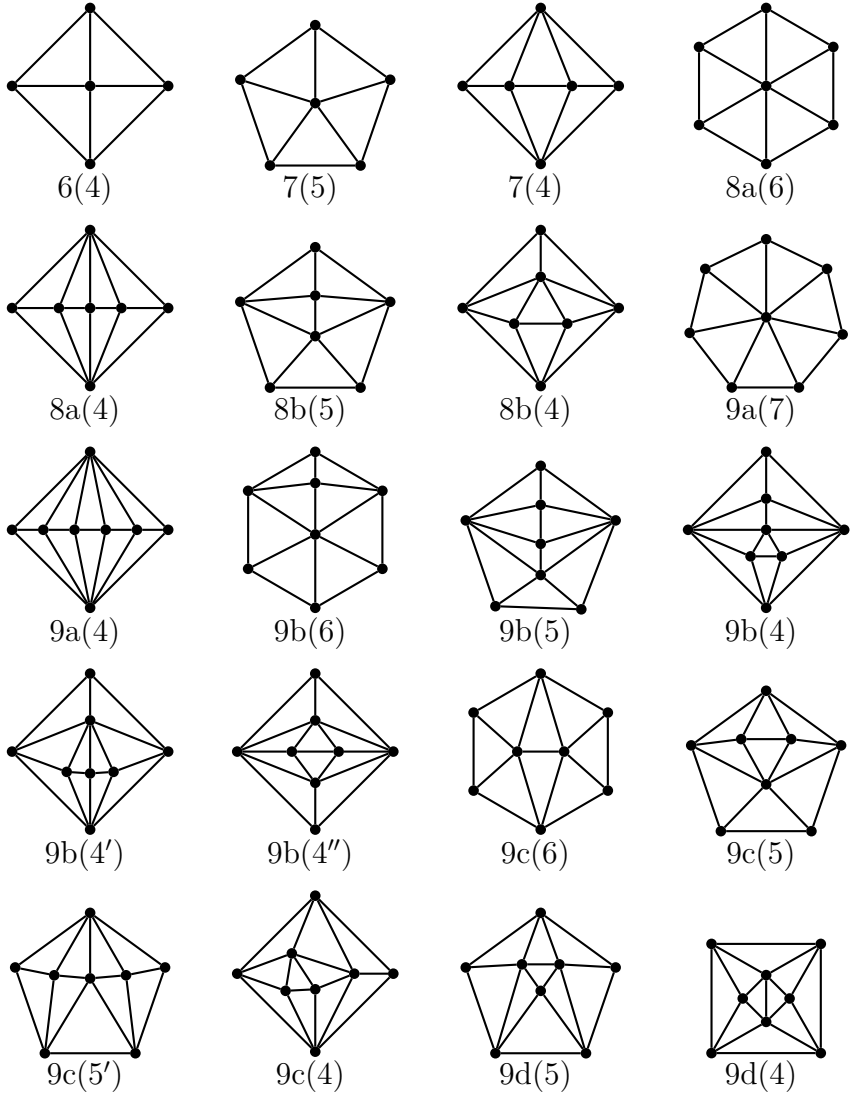
\begin{figure}
\centering
\psset{unit=0.0013\textwidth}

\parbox[b]{.15\textwidth}{%
\centering
\begin{pspicture}(-51,-51)(51,51)
\cnode*(0,50){2pt}{e1}
\cnode*(50,0){2pt}{d1}
\cnode*(0,-50){2pt}{c1}
\cnode*(-50,0){2pt}{b1}
\cnode*(0,0){2pt}{f1}
\ncline{e1}{d1}
\ncline{e1}{b1}
\ncline{e1}{f1}
\ncline{d1}{c1}
\ncline{d1}{f1}
\ncline{c1}{b1}
\ncline{c1}{f1}
\ncline{b1}{f1}
\rput(0,-66){\rnode{A}{6(4)
}}
\end{pspicture}}
\hspace{.02\textwidth}
\parbox[b]{.15\textwidth}{%
\centering
\begin{pspicture}(-48,-40)(48,50)
\cnode*(0,50){2pt}{d1}
\cnode*(48,15){2pt}{g1}
\cnode*(29,-40){2pt}{f1}
\cnode*(-29,-40){2pt}{b1}
\cnode*(-48,15){2pt}{a1}
\cnode*(0,0){2pt}{e1}
\ncline{d1}{g1}
\ncline{d1}{a1}
\ncline{d1}{e1}
\ncline{g1}{f1}
\ncline{g1}{e1}
\ncline{f1}{b1}
\ncline{f1}{e1}
\ncline{b1}{a1}
\ncline{b1}{e1}
\ncline{a1}{e1}
\rput(0,-55){\rnode{A}{7(5)
}}
\end{pspicture}}
\hspace{.02\textwidth}
\parbox[b]{.15\textwidth}{%
\centering
\begin{pspicture}(-51,-51)(51,51)
\cnode*(0,50){2pt}{e1}
\cnode*(50,0){2pt}{d1}
\cnode*(0,-50){2pt}{c1}
\cnode*(-50,0){2pt}{b1}
\cnode*(-20,0){2pt}{f1}
\cnode*(20,0){2pt}{g1}
\ncline{e1}{d1}
\ncline{e1}{b1}
\ncline{e1}{f1}
\ncline{e1}{g1}
\ncline{d1}{c1}
\ncline{d1}{g1}
\ncline{c1}{b1}
\ncline{c1}{f1}
\ncline{c1}{g1}
\ncline{b1}{f1}
\ncline{f1}{g1}
\rput(0,-66){\rnode{A}{7(4)
}}
\end{pspicture}}
\hspace{.02\textwidth}
\parbox[b]{.15\textwidth}{%
\centering
\begin{pspicture}(-43,-51)(43,51)
\cnode*(0,50){2pt}{d1}
\cnode*(43,25){2pt}{h1}
\cnode*(43,-25){2pt}{g1}
\cnode*(0,-50){2pt}{f1}
\cnode*(-43,-25){2pt}{b1}
\cnode*(-43,25){2pt}{a1}
\cnode*(0,0){2pt}{e1}
\ncline{d1}{h1}
\ncline{d1}{a1}
\ncline{d1}{e1}
\ncline{h1}{g1}
\ncline{h1}{e1}
\ncline{g1}{f1}
\ncline{g1}{e1}
\ncline{f1}{b1}
\ncline{f1}{e1}
\ncline{b1}{a1}
\ncline{b1}{e1}
\ncline{a1}{e1}
\rput(0,-66){\rnode{A}{8a(6)
}}
\end{pspicture}}

\vspace{.05\textwidth}

\parbox[b]{.15\textwidth}{%
\centering
\begin{pspicture}(-51,-51)(51,51)
\cnode*(0,50){2pt}{e1}
\cnode*(50,0){2pt}{d1}
\cnode*(0,-50){2pt}{c1}
\cnode*(-50,0){2pt}{b1}
\cnode*(-20,0){2pt}{f1}
\cnode*(0,0){2pt}{g1}
\cnode*(20,0){2pt}{h1}
\ncline{e1}{d1}
\ncline{e1}{b1}
\ncline{e1}{f1}
\ncline{e1}{g1}
\ncline{e1}{h1}
\ncline{d1}{c1}
\ncline{d1}{h1}
\ncline{c1}{b1}
\ncline{c1}{f1}
\ncline{c1}{g1}
\ncline{c1}{h1}
\ncline{b1}{f1}
\ncline{f1}{g1}
\ncline{g1}{h1}
\rput(0,-66){\rnode{A}{8a(4)
}}
\end{pspicture}}
\hspace{.02\textwidth}
\parbox[b]{.15\textwidth}{%
\centering
\begin{pspicture}(-48,-40)(48,50)
\cnode*(48,15){2pt}{c1}
\cnode*(29,-40){2pt}{g1}
\cnode*(-29,-40){2pt}{f1}
\cnode*(-48,15){2pt}{e1}
\cnode*(0,50){2pt}{a1}
\cnode*(0,19){2pt}{d1}
\cnode*(0,-7){2pt}{h1}
\ncline{c1}{g1}
\ncline{c1}{a1}
\ncline{c1}{d1}
\ncline{c1}{h1}
\ncline{g1}{f1}
\ncline{g1}{h1}
\ncline{f1}{e1}
\ncline{f1}{h1}
\ncline{e1}{a1}
\ncline{e1}{d1}
\ncline{e1}{h1}
\ncline{a1}{d1}
\ncline{d1}{h1}
\rput(0,-55){\rnode{A}{8b(5)
}}
\end{pspicture}}
\hspace{.02\textwidth}
\parbox[b]{.15\textwidth}{%
\centering
\begin{pspicture}(-51,-51)(51,51)
\cnode*(-50,0){2pt}{e1}
\cnode*(0,50){2pt}{d1}
\cnode*(50,0){2pt}{c1}
\cnode*(0,-50){2pt}{b1}
\cnode*(-17,-10){2pt}{f1}
\cnode*(17,-10){2pt}{g1}
\cnode*(0,20){2pt}{h1}
\ncline{e1}{d1}
\ncline{e1}{b1}
\ncline{e1}{f1}
\ncline{e1}{h1}
\ncline{d1}{c1}
\ncline{d1}{h1}
\ncline{c1}{b1}
\ncline{c1}{g1}
\ncline{c1}{h1}
\ncline{b1}{f1}
\ncline{b1}{g1}
\ncline{f1}{g1}
\ncline{f1}{h1}
\ncline{g1}{h1}
\rput(0,-66){\rnode{A}{8b(4)
}}
\end{pspicture}}
\hspace{.02\textwidth}
\parbox[b]{.15\textwidth}{%
\centering
\begin{pspicture}(-49,-45)(49,50)
\cnode*(0,50){2pt}{d1}
\cnode*(39,31){2pt}{i1}
\cnode*(49,-11){2pt}{h1}
\cnode*(21,-45){2pt}{g1}
\cnode*(-22,-45){2pt}{f1}
\cnode*(-49,-10){2pt}{b1}
\cnode*(-39,31){2pt}{a1}
\cnode*(0,0){2pt}{e1}
\ncline{d1}{i1}
\ncline{d1}{a1}
\ncline{d1}{e1}
\ncline{i1}{h1}
\ncline{i1}{e1}
\ncline{h1}{g1}
\ncline{h1}{e1}
\ncline{g1}{f1}
\ncline{g1}{e1}
\ncline{f1}{b1}
\ncline{f1}{e1}
\ncline{b1}{a1}
\ncline{b1}{e1}
\ncline{a1}{e1}
\rput(0,-60){\rnode{A}{9a(7)
}}
\end{pspicture}}

\vspace{.05\textwidth}

\parbox[b]{.15\textwidth}{%
\centering
\begin{pspicture}(-51,-51)(51,51)
\cnode*(0,50){2pt}{e1}
\cnode*(50,0){2pt}{d1}
\cnode*(0,-50){2pt}{c1}
\cnode*(-50,0){2pt}{b1}
\cnode*(-30,0){2pt}{f1}
\cnode*(-10,0){2pt}{g1}
\cnode*(10,0){2pt}{h1}
\cnode*(30,0){2pt}{i1}
\ncline{e1}{d1}
\ncline{e1}{b1}
\ncline{e1}{f1}
\ncline{e1}{g1}
\ncline{e1}{h1}
\ncline{e1}{i1}
\ncline{d1}{c1}
\ncline{d1}{i1}
\ncline{c1}{b1}
\ncline{c1}{f1}
\ncline{c1}{g1}
\ncline{c1}{h1}
\ncline{c1}{i1}
\ncline{b1}{f1}
\ncline{f1}{g1}
\ncline{g1}{h1}
\ncline{h1}{i1}
\rput(0,-66){\rnode{A}{9a(4)
}}
\end{pspicture}}
\hspace{.02\textwidth}
\parbox[b]{.15\textwidth}{%
\centering
\begin{pspicture}(-43,-51)(43,51)
\cnode*(-43,-25){2pt}{d1}
\cnode*(-43,25){2pt}{h1}
\cnode*(0,50){2pt}{g1}
\cnode*(43,25){2pt}{f1}
\cnode*(43,-25){2pt}{b1}
\cnode*(0,-50){2pt}{a1}
\cnode*(0,-3){2pt}{e1}
\cnode*(0,30){2pt}{i1}
\ncline{d1}{h1}
\ncline{d1}{a1}
\ncline{d1}{e1}
\ncline{h1}{g1}
\ncline{h1}{e1}
\ncline{h1}{i1}
\ncline{g1}{f1}
\ncline{g1}{i1}
\ncline{f1}{b1}
\ncline{f1}{e1}
\ncline{f1}{i1}
\ncline{b1}{a1}
\ncline{b1}{e1}
\ncline{a1}{e1}
\ncline{e1}{i1}
\rput(0,-66){\rnode{A}{9b(6)
}}
\end{pspicture}}
\hspace{.02\textwidth}
\parbox[b]{.15\textwidth}{%
\centering
\begin{pspicture}(-48,-42)(48,50)
\cnode*(48,15){2pt}{c1}
\cnode*(26,-42){2pt}{g1}
\cnode*(-29,-40){2pt}{i1}
\cnode*(-48,15){2pt}{e1}
\cnode*(0,50){2pt}{b1}
\cnode*(0,25){2pt}{a1}
\cnode*(0,0){2pt}{d1}
\cnode*(0,-20){2pt}{h1}
\ncline{c1}{g1}
\ncline{c1}{b1}
\ncline{c1}{a1}
\ncline{c1}{d1}
\ncline{c1}{h1}
\ncline{g1}{i1}
\ncline{g1}{h1}
\ncline{i1}{e1}
\ncline{i1}{h1}
\ncline{e1}{b1}
\ncline{e1}{a1}
\ncline{e1}{d1}
\ncline{e1}{h1}
\ncline{b1}{a1}
\ncline{a1}{d1}
\ncline{d1}{h1}
\rput(0,-57){\rnode{A}{9b(5)
}}
\end{pspicture}}
\hspace{.02\textwidth}
\parbox[b]{.15\textwidth}{%
\centering
\begin{pspicture}(-51,-51)(51,51)
\cnode*(50,0){2pt}{c1}
\cnode*(0,-50){2pt}{f1}
\cnode*(-50,0){2pt}{e1}
\cnode*(0,50){2pt}{a1}
\cnode*(0,20){2pt}{d1}
\cnode*(10,-17){2pt}{g1}
\cnode*(0,0){2pt}{h1}
\cnode*(-10,-17){2pt}{i1}
\ncline{c1}{f1}
\ncline{c1}{a1}
\ncline{c1}{d1}
\ncline{c1}{g1}
\ncline{c1}{h1}
\ncline{f1}{e1}
\ncline{f1}{g1}
\ncline{f1}{i1}
\ncline{e1}{a1}
\ncline{e1}{d1}
\ncline{e1}{h1}
\ncline{e1}{i1}
\ncline{a1}{d1}
\ncline{d1}{h1}
\ncline{g1}{h1}
\ncline{g1}{i1}
\ncline{h1}{i1}
\rput(0,-66){\rnode{A}{9b(4)
}}
\end{pspicture}}

\vspace{.05\textwidth}

\parbox[b]{.15\textwidth}{%
\centering
\begin{pspicture}(-51,-51)(51,51)
\cnode*(-50,0){2pt}{h1}
\cnode*(0,50){2pt}{i1}
\cnode*(50,0){2pt}{f1}
\cnode*(0,-50){2pt}{c1}
\cnode*(0,-14){2pt}{a1}
\cnode*(15,-13){2pt}{b1}
\cnode*(-15,-13){2pt}{d1}
\cnode*(0,20){2pt}{e1}
\ncline{h1}{i1}
\ncline{h1}{c1}
\ncline{h1}{d1}
\ncline{h1}{e1}
\ncline{i1}{f1}
\ncline{i1}{e1}
\ncline{f1}{c1}
\ncline{f1}{b1}
\ncline{f1}{e1}
\ncline{c1}{a1}
\ncline{c1}{b1}
\ncline{c1}{d1}
\ncline{a1}{b1}
\ncline{a1}{d1}
\ncline{a1}{e1}
\ncline{b1}{e1}
\ncline{d1}{e1}
\rput(0,-66){\rnode{A}{9b$(4')$
}}
\end{pspicture}}
\hspace{.02\textwidth}
\parbox[b]{.15\textwidth}{%
\centering
\begin{pspicture}(-51,-51)(51,51)
\cnode*(-50,0){2pt}{e1}
\cnode*(0,50){2pt}{d1}
\cnode*(50,0){2pt}{c1}
\cnode*(0,-50){2pt}{b1}
\cnode*(0,-20){2pt}{f1}
\cnode*(15,0){2pt}{g1}
\cnode*(0,20){2pt}{h1}
\cnode*(-15,0){2pt}{i1}
\ncline{e1}{d1}
\ncline{e1}{b1}
\ncline{e1}{f1}
\ncline{e1}{h1}
\ncline{e1}{i1}
\ncline{d1}{c1}
\ncline{d1}{h1}
\ncline{c1}{b1}
\ncline{c1}{f1}
\ncline{c1}{g1}
\ncline{c1}{h1}
\ncline{b1}{f1}
\ncline{f1}{g1}
\ncline{f1}{i1}
\ncline{g1}{h1}
\ncline{g1}{i1}
\ncline{h1}{i1}
\rput(0,-66){\rnode{A}{9b$(4'')$
}}
\end{pspicture}}
\hspace{.02\textwidth}
\parbox[b]{.15\textwidth}{%
\centering
\begin{pspicture}(-43,-51)(43,51)
\cnode*(43,-25){2pt}{d1}
\cnode*(0,-50){2pt}{i1}
\cnode*(-43,-25){2pt}{h1}
\cnode*(-43,25){2pt}{g1}
\cnode*(0,50){2pt}{b1}
\cnode*(43,25){2pt}{a1}
\cnode*(15,0){2pt}{e1}
\cnode*(-15,0){2pt}{f1}
\ncline{d1}{i1}
\ncline{d1}{a1}
\ncline{d1}{e1}
\ncline{i1}{h1}
\ncline{i1}{e1}
\ncline{i1}{f1}
\ncline{h1}{g1}
\ncline{h1}{f1}
\ncline{g1}{b1}
\ncline{g1}{f1}
\ncline{b1}{a1}
\ncline{b1}{e1}
\ncline{b1}{f1}
\ncline{a1}{e1}
\ncline{e1}{f1}
\rput(0,-66){\rnode{A}{9c(6)
}}
\end{pspicture}}
\hspace{.02\textwidth}
\parbox[b]{.15\textwidth}{%
\centering
\begin{pspicture}(-48,-40)(48,50)
\cnode*(-48,15){2pt}{b1}
\cnode*(0,50){2pt}{f1}
\cnode*(48,15){2pt}{i1}
\cnode*(29,-40){2pt}{d1}
\cnode*(-29,-40){2pt}{a1}
\cnode*(0,-10){2pt}{c1}
\cnode*(-16,19){2pt}{g1}
\cnode*(16,19){2pt}{h1}
\ncline{b1}{f1}
\ncline{b1}{a1}
\ncline{b1}{c1}
\ncline{b1}{g1}
\ncline{f1}{i1}
\ncline{f1}{g1}
\ncline{f1}{h1}
\ncline{i1}{d1}
\ncline{i1}{c1}
\ncline{i1}{h1}
\ncline{d1}{a1}
\ncline{d1}{c1}
\ncline{a1}{c1}
\ncline{c1}{g1}
\ncline{c1}{h1}
\ncline{g1}{h1}
\rput(0,-55){\rnode{A}{9c(5)
}}
\end{pspicture}}

\vspace{.05\textwidth}

\parbox[b]{.15\textwidth}{%
\centering
\begin{pspicture}(-48,-40)(48,50)
\cnode*(0,50){2pt}{c1}
\cnode*(48,15){2pt}{g1}
\cnode*(29,-40){2pt}{f1}
\cnode*(-29,-40){2pt}{e1}
\cnode*(-48,15){2pt}{a1}
\cnode*(-23,10){2pt}{d1}
\cnode*(23,10){2pt}{h1}
\cnode*(0,8){2pt}{i1}
\ncline{c1}{g1}
\ncline{c1}{a1}
\ncline{c1}{d1}
\ncline{c1}{h1}
\ncline{c1}{i1}
\ncline{g1}{f1}
\ncline{g1}{h1}
\ncline{f1}{e1}
\ncline{f1}{h1}
\ncline{f1}{i1}
\ncline{e1}{a1}
\ncline{e1}{d1}
\ncline{e1}{i1}
\ncline{a1}{d1}
\ncline{d1}{i1}
\ncline{h1}{i1}
\rput(0,-55){\rnode{A}{9c$(5')$
}}
\end{pspicture}}
\hspace{.02\textwidth}
\parbox[b]{.15\textwidth}{%
\centering
\begin{pspicture}(-51,-51)(51,51)
\cnode*(0,50){2pt}{e1}
\cnode*(50,0){2pt}{d1}
\cnode*(0,-50){2pt}{c1}
\cnode*(-50,0){2pt}{b1}
\cnode*(-15,13){2pt}{f1}
\cnode*(-19,-11){2pt}{g1}
\cnode*(0,-10){2pt}{h1}
\cnode*(25,0){2pt}{i1}
\ncline{e1}{d1}
\ncline{e1}{b1}
\ncline{e1}{f1}
\ncline{e1}{i1}
\ncline{d1}{c1}
\ncline{d1}{i1}
\ncline{c1}{b1}
\ncline{c1}{g1}
\ncline{c1}{h1}
\ncline{c1}{i1}
\ncline{b1}{f1}
\ncline{b1}{g1}
\ncline{f1}{g1}
\ncline{f1}{h1}
\ncline{f1}{i1}
\ncline{g1}{h1}
\ncline{h1}{i1}
\rput(0,-66){\rnode{A}{9c(4)
}}
\end{pspicture}}
\hspace{.02\textwidth}
\parbox[b]{.15\textwidth}{%
\centering
\begin{pspicture}(-48,-40)(48,50)
\cnode*(-48,15){2pt}{f1}
\cnode*(0,50){2pt}{e1}
\cnode*(48,15){2pt}{d1}
\cnode*(29,-40){2pt}{c1}
\cnode*(-29,-40){2pt}{b1}
\cnode*(-12,17){2pt}{g1}
\cnode*(0,0){2pt}{h1}
\cnode*(12,17){2pt}{i1}
\ncline{f1}{e1}
\ncline{f1}{b1}
\ncline{f1}{g1}
\ncline{e1}{d1}
\ncline{e1}{g1}
\ncline{e1}{i1}
\ncline{d1}{c1}
\ncline{d1}{i1}
\ncline{c1}{b1}
\ncline{c1}{h1}
\ncline{c1}{i1}
\ncline{b1}{g1}
\ncline{b1}{h1}
\ncline{g1}{h1}
\ncline{g1}{i1}
\ncline{h1}{i1}
\rput(0,-55){\rnode{A}{9d(5)
}}
\end{pspicture}}
\hspace{.02\textwidth}
\parbox[b]{.15\textwidth}{%
\centering
\begin{pspicture}(-35,-35)(35,35)
\cnode*(-35,35){2pt}{e1}
\cnode*(35,35){2pt}{i1}
\cnode*(35,-35){2pt}{c1}
\cnode*(-35,-35){2pt}{a1}
\cnode*(0,-15){2pt}{b1}
\cnode*(-15,0){2pt}{f1}
\cnode*(0,15){2pt}{g1}
\cnode*(15,0){2pt}{h1}
\ncline{e1}{i1}
\ncline{e1}{a1}
\ncline{e1}{f1}
\ncline{e1}{g1}
\ncline{i1}{c1}
\ncline{i1}{g1}
\ncline{i1}{h1}
\ncline{c1}{a1}
\ncline{c1}{b1}
\ncline{c1}{h1}
\ncline{a1}{b1}
\ncline{a1}{f1}
\ncline{b1}{f1}
\ncline{b1}{g1}
\ncline{b1}{h1}
\ncline{f1}{g1}
\ncline{g1}{h1}
\rput(0,-50){\rnode{A}{9d(4)
}}
\end{pspicture}}

\vspace{.05\textwidth}

\caption{$L_vu$ when the degree of $u = 6,7,8,9$}
\label{fig:WalkupFig3}
\end{figure}

\begin{thm}  If $K$ is $g_2$-irreducible, $u$ a vertex of degree less than $10$ and $v$ a vertex in the link of $u$, 
then the one-skeleton of the link of $u$ with $v$ and its incident edges removed is exactly one of those 
in Figures~\ref{fig:WalkupFig3}-6(4)--\ref{fig:WalkupFig3}-9d(4).
\end{thm}

From here on we write ``$L_v u$ is of type" to mean that $v$ is in the link of $u$ in a $g_2$-irreducible 
triangulation, and the one-skeleton of the link of $u$ with $v$ and its incident edges removed is the 
referenced figure.

\begin{thm}  \label{thm:minvertexdeg} Let $K$ be a $g_2$-irreducible triangulation.  Then there exists a triangulation $K^\prime$ 
which is homeomorphic to $K$, has the same $f$-vector as $K,$ and whose links satisfy the following:
\begin{itemize}
\addtolength{\itemsep}{-.5mm}
 \item If $L_v u$ is of type 6(4), then\, ${\deg}(v)\geq 10$.
 \item If $L_v u$ is of type 7(5), then\, ${\deg}(v)\geq 12$.
 \item If $L_v u$ is of type 8a(6), then\, ${\deg}(v)\geq 14$.
 \item If $L_v u$ is of type 8b(5), then\, ${\deg}(v)\geq 11$.
 \item If all of the vertices of $K^\prime$ have degree at least $9,$ 
       then there exists at least two vertices of degree at least $10$, 
       or there exists at least one vertex of degree at least $11.$
\addtolength{\itemsep}{.5mm}
\end{itemize}
\end{thm}

\textbf{Proof of Theorem \ref{thm:newwalkup}:}  
Let $K$ be $g_2$-irreducible.  We can assume that $K$ satisfies the conclusions of the previous theorem.  Let $(u,v)$ be an ordered pair of vertices of $K$ 
which form an edge.  Define $\lambda(u,v)$ as follows:

\begin{itemize}
\addtolength{\itemsep}{-.5mm}
 \item $\lambda(u,v) = \frac{3}{4}$\, if $L_v u$ is of type $6(4).$
 \item $\lambda(u,v) = 1$\, if $L_v u$ is of type $7(5).$
 \item $\lambda(u,v) = \frac{3}{4}$\, if $L_v u$ is of type $8a(6).$
 \item $\lambda(u,v) = \frac{5}{8}$\, if $L_v u$ is of type $8b(5).$
 \item $\lambda(u,v) = \frac{1}{2}$\, if $L_v u$ is of type $7(4), 8a(4), 8b(4),$ or if the degree of $u$ is $9$.
 \item $\lambda(u,v) = 1- \lambda(v,u)$\, if the degree of $u$ is at least $10$ and the degree of $v$ is $9$ or less. 
 \item $\lambda(u,v) = \frac{1}{2}$\, otherwise.
\addtolength{\itemsep}{.5mm}
\end{itemize}

Define 
$$\mu(u) = \sum_{v\in {\rm Lk}\,u} \lambda(u,v) - \tfrac{9}{2}.$$

\noindent
By construction,
$$ \sum_{u \in K} \mu(u) = f_1(K) - \tfrac{9}{2} f_0(K).$$

Suppose that $u$ is a vertex of degree $m.$   If $v$ is in the link of $u$ and $L_u v$ is of type $6(4)$, $7(5)$, $8a(6)$, or $8b(5)$, 
then Theorem \ref{thm:minvertexdeg} implies that the degree of $u$ is at least ten. Therefore, in the link of $u$  each triangle has at most one vertex $v$ such that $L_u v$ is one of these four types.  
Let $n_{6(4)}, n_{7(4)}, n_{7(5)},$ etc., be the number of vertices $v$ in the link of $u$ of such that $L_u v$ is of type $6(4), 7(4), 7(5),$ etc.  
Since the link of $u$ has $2m-4$ triangles, the link of $u$ must satisfy the integer constraint 
\begin{equation} \label{integerconstraint}
4 n_{6(4)} + 5 n_{7(5)} + 6 n_{8a(6)} + 5 n_{8b(5)}  \le  2m-4.
\end{equation}

\noindent  
The minimum potential value of $\mu(u)$ is the minimum of 
\begin{equation} \label{mu(u)}
-\frac{1}{4} n_{6(4)} - \frac{1}{2} n_{7(5)} - \frac{1}{4} n_{8a(6)} - \frac{1}{8} n_{8b(5)} + \frac{m-9}{2}
\end{equation}
under the above constraint and a few others discussed below.  Now we determine lower bounds for $\mu(u)$ for a variety of values of $m.$ 

\begin{itemize}
  \item $m < 10.$  Then by definition $\mu(u) =0.$  
  \item $m=10.$  $\mu(u) \ge \frac{1}{4}$\, \cite[Lemma 11.9]{Walkup1970}.
  \item $m=11.$ $\mu(u) \ge \frac{1}{2}$\, \cite[Lemma 11.9]{Walkup1970}.
  \item $m=12.$  $\mu(u) \ge \frac{1}{2}$\, \cite[Lemma 11.9]{Walkup1970}.
  \item $m=13.$ Walkup proved that $n_{7(5)} \le 3$ in this case.  With this additional restriction, 
        the minimal value of (\ref{mu(u)}) subject to the integral constraint (\ref{integerconstraint}) is $\frac{1}{4}.$  
        This value occurs when $n_{6(4)}=3$ and $n_{7(5)}=2,$ or $n_{6(4)}=1$ and $n_{7(5)} = 3.$  
        In all other cases $\mu(u) > \frac{1}{4}.$
  \item $m=14.$  The minimal  value of (\ref{mu(u)}) subject to the integral constraint (\ref{integerconstraint}) is  $\frac{1}{4}$ 
        and this only occurs if $n_{6(4)} = 1$ and $n_{7(5)}=4.$  In all other cases, $\mu(u) > \frac{1}{4}.$
  \item $m \ge 15.$  Even without integer considerations, $\mu(u) > \frac{1}{4}.$
\end{itemize}

For notational purposes, define 
$$\mu(K) =  \sum_{u \in K} \mu(u) = f_1(K) - \tfrac{9}{2} f_0(K).$$
From above we know that $\mu(u) \ge 0 $ for all vertices $u.$  If $K$ has no vertices of degree less than $9,$ 
then by the last line of Theorem \ref{thm:minvertexdeg} there exists at least two vertices which contribute at least $1/2$ 
to $\mu(K)$ or one vertex which contributes at least $1,$ so $\mu(K) \ge 1.$  So suppose $K$ has a vertex of degree less than nine.  There are four possibilities.

\begin{enumerate}
\addtolength{\itemsep}{-.5mm}
 \item $K$ has a vertex of degree six. Then the six vertices of the link of this vertex all contribute at least $1/4$ to $\mu(K).$
 \item $K$ has a vertex of degree seven.  Consider the two vertices of type $7(5)$ whose existence is now guaranteed.  
       Each of these either adds more than $1/4$ to $\mu(K)$ or imply the existence of a vertex of degree six.
 \item $K$ has a vertex whose link is of type $8a.$  The same argument as the case of a vertex of degree seven applies.
 \item $K$ has a vertex whose link is of type $8b.$  Then $K$ has at least four vertices of type $8b(5)$ 
       each of which either satisfy $\mu(u) > 1/4$ or imply the existence of a vertex of degree six. \hfill$\Box$ 
\addtolength{\itemsep}{.5mm}
\end{enumerate}

\section{Enumeration of \mathversion{bold}$g_2$-Irreducible\mathversion{normal} Triangulations}
\label{sec:enumeration}

A $3$-manifold $M$ is \emph{irreducible} if every embedded $2$-sphere in $M$ bounds a $3$-ball in $M$. 
In particular, if a triangulation $K$ of an irreducible $3$-manifold $M$ has a missing facet, then
$K=K_1\# K_2$ with one part homeomorphic to $M$ and the other homeomorphic to $S^3$.
As a consequence, every $g_2$-minimal triangulation $K$ of an irreducible $3$-manifold $M$, 
different from $S^3$, is either $g_2$-irreducible or is obtained from a $g_2$-irreducible 
triangulation $K^\prime$ of $M$ by successive stacking operations $S$.

As already mentioned in the previous section, for every fixed $\Gamma$ there are only finitely many $3$-manifolds $M$
such that $\Gamma(M)\leq\Gamma$ \cite{Swartz2008}: If $M$ is a $3$-manifold with $\Gamma(M)\leq\Gamma$, 
then $M$ has a triangulation $K$ with $f$-vector $f=(f_0,f_1,f_2,f_3)$ such that $f_1-4f_0+10\leq\Gamma$.
If $K$ is $g_2$-irreducible, then the additional restriction~(\ref{eq:newwalkup}) holds, $f_1>\tfrac{9}{2} f_0+\tfrac{1}{2}$. 
These two inequalities together with the trivial inequality $f_1\leq\binom{f_0}{2}$ allow for only 
finitely many tuples $(f_0,f_1)$. Hence, there are only finitely many $g_2$-irreducible 
triangulations $K$ with $g_2(K)\leq\Gamma$. This directly implies
that there are only finitely many irreducible $3$-manifolds $M$ with $\Gamma(M)\leq\Gamma$.
If $K$ is a $g_2$-minimal triangulation of a non-irreducible $3$-manifold $M$ with $\Gamma(M)\leq\Gamma$,
then $K$ is either $g_2$-irreducible, in which case there are only finitely many such triangulations,
or $K$ is of the form $K_1 \# K_2$ or $HK^\prime$, where the component triangulations realize their minimum $g_2$.
These components are either $g_2$-irreducible or can further be split up or reduced by deleting a handle. 
Since $g_2(K_1\# K_2) = g_2(K_1)+g_2(K_2)$ and $g_2(HK^\prime)=g_2(K^\prime)+10$,
it follows that there are at most finitely many non-irreducible $3$-manifolds $M$ with $\Gamma(M)\leq\Gamma$
and therefore only finitely many $3$-man\-i\-folds $M$ with $\Gamma(M)\leq\Gamma$

\begin{figure}
\begin{center}
\psfrag{1}{1}
\psfrag{2}{2}
\psfrag{3}{3}
\psfrag{4}{4}
\psfrag{5}{5}
\psfrag{6}{6}
\psfrag{7}{7}
\psfrag{8}{8}
\psfrag{9}{9}
\psfrag{10}{10}
\psfrag{11}{11}
\psfrag{12}{12}
\psfrag{13}{13}
\psfrag{14}{14}
\psfrag{15}{15}
\psfrag{16}{16}
\psfrag{17}{17}
\psfrag{18}{18}
\psfrag{19}{19}
\psfrag{20}{20}
\psfrag{30}{30}
\psfrag{40}{40}
\psfrag{50}{50}
\psfrag{60}{60}
\psfrag{70}{70}
\psfrag{80}{80}
\psfrag{90}{90}
\psfrag{100}{100}
\psfrag{110}{110}
\psfrag{-10}{$-10$}
\psfrag{(5,10)}{$(5,10)$}
\psfrag{(6,15)}{$(6,15)$}
\psfrag{(7,21)}{$(7,21)$}
\psfrag{(8,28)}{$(8,28)$}
\psfrag{(9,36)}{$(9,36)$}
\psfrag{(10,45)}{$(10,45)$}
\psfrag{(11,55)}{$(11,55)$}
\psfrag{(12,66)}{$(12,66)$}
\psfrag{(13,78)}{$(13,78)$}
\psfrag{(14,91)}{$(14,91)$}
\psfrag{(15,105)}{$(15,105)$}
\psfrag{f0}{$f_0$}
\psfrag{f1}{$f_1$}
\psfrag{f1=4.5f0+0.5}{$f_1=4.5f_0+0.5$}
\psfrag{f1=4f0+10}{$f_1=4f_0+10$\, ($g_2=20$)}
\psfrag{f1=4f0}{$f_1=4f_0$\, ($g_2=10$)}
\psfrag{f1=4f0-10}{$f_1=4f_0-10$\, ($g_2=0$)}
\psfrag{binom}{$f_1=\binom{f_0}{2}$}
\includegraphics[width=.899\linewidth]{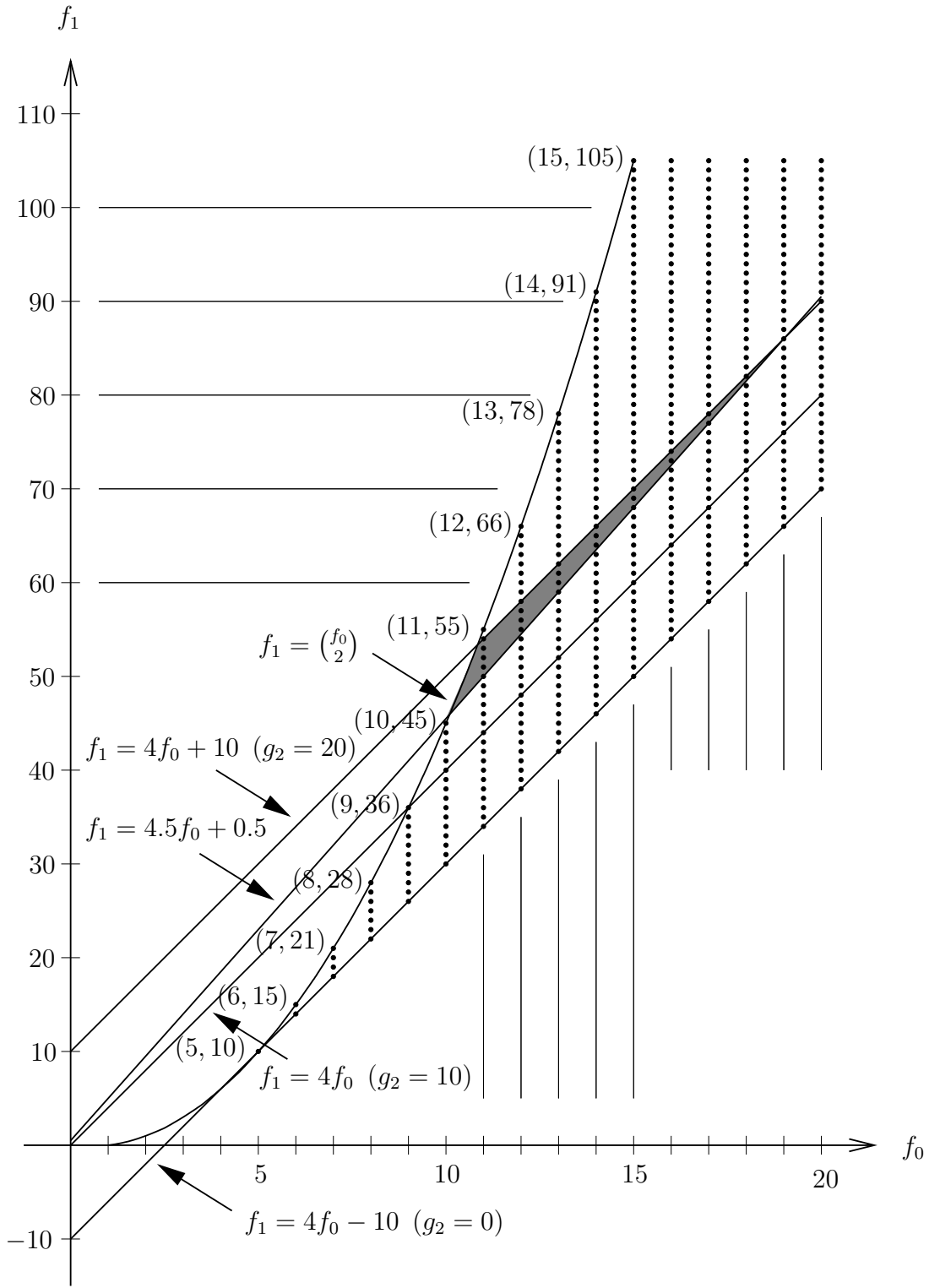}
\end{center}
\caption{Range (in grey) for $g_2$-irreducible triangulations of $3$-manifolds with $g_2\leq 20$.
The range is bounded by the inequalities $f_1>\tfrac{9}{2} f_0+\tfrac{1}{2}$,
$f_1-4f_0+10\leq 20$, and $f_1\leq\binom{f_0}{2}$.}
\label{fig:diagram}
\end{figure}

Figure~\ref{fig:diagram} displays in grey the admissible range for tuples $(f_0,f_1)$
that can occur for $g_2$-irreducible triangulations of $3$-manifolds with $\Gamma\leq 20$. 
These are precisely the tuples $(11,51)$, $(11,52)$, $(11,53)$, $(11,54)$,
$(12,55)$, $(12,56)$, $(12,57)$, $(12,58)$, $(13,60)$, $(13,61)$, $(13,62)$, $(14,64)$, $(14,65)$,
$(14,66)$, $(15,69)$, $(15,70)$, $(16,73)$, $(16,74)$, $(17,78)$, and $(18,82)$.

We conducted exhaustive computer searches to find all the $g_2$-irreducible triangulations of $3$-manifolds 
with $g_2\leq 20$ and candidates for all $g_2$-irreducible triangulations of $3$-manifolds with $f_0 \leq 15$.

The $3$-manifolds were constructed using the lexicographic enumeration technique described in \cite{SulankeLutz2006pre}.
This technique constructs $3$-manifolds one facet at a time which allows local properties to be tested 
before the $3$-manifolds are completely constructed. By checking local properties provided by Walkup \cite{Walkup1970} 
the searches can be pruned sufficiently to make them feasible.

\begin{thm}  {\rm (Walkup \cite[10.1]{Walkup1970})} \label{thm:lemma101} 
             Let $K$ be a $g_2$-irreducible triangulation and $(u,v)$ be an edge of $K$.  
             Then ${\rm Lk}\,(u,v)$ contains at least $4$ vertices.
\end{thm}

\begin{thm}  {\rm (Walkup \cite[10.2]{Walkup1970})} \label{thm:lemma102} 
             Let $K$ be a $g_2$-irreducible triangulation and $(u,v)$ be an edge of $K$.  
             Then ${\rm Lk}\,u \cap {\rm Lk}\,v - {\rm Lk}\,(u,v)$ is nonempty.
\end{thm}

\begin{thm}  {\rm (Walkup \cite[10.4]{Walkup1970})} \label{thm:lemma104} 
             Let $K$ be a $g_2$-irreducible triangulation and $u$ be a vertex of $K$.  
             Suppose ${\rm Lk}\,u$ contains the boundary complex of a $2$-simplex $(a,b,c)$ as a subcomplex.  
             Then  ${\rm Lk}\,u$ must also contain the $2$-simplex  $(a,b,c)$.
\end{thm}

\begin{thm}  {\rm (Walkup \cite[11.1]{Walkup1970})} \label{thm:lemma111} 
             Let $K$ be a $g_2$-irreducible triangulation and  $(u,v)$ be an edge of $K$.  
             Suppose ${\rm Lk}\,u \cap {\rm Lk}\,v - {\rm Lk}\,(u,v) = \{w\}$.  
             Then ${\rm Lk}\,(u,w)$ contains at least as many vertices as ${\rm Lk}\,(u,v)$.
\end{thm}

To find all the candidates for $g_2$-irreducible triangulations of $3$-manifolds with $f_0 \leq 15$
Theorems \ref{thm:lemma101}, \ref{thm:lemma102}, \ref{thm:lemma104}, \ref{thm:lemma111}
were used to prune the searches.  
Run times for $11,12,13,14,15$ vertices were 2 seconds, 100 seconds, 3 hours, 60 days, and 7 years, respectively.
Results from these runs are given in Table~\ref{tbl:g2minimal_examples}.

\begin{table}
\small\centering
\defaultaddspace=0.125em
\caption{Triangulations found with $5\leq f_0\leq 15$ vertices after
implementing Walkup's Lemmas 10.1, 10.2, 10.4, and 11.1.}\label{tbl:g2minimal_examples}
\begin{tabular}{@{}l@{\hspace{12mm}}r@{\hspace{12mm}}l@{\hspace{10mm}}l@{\hspace{10mm}}r@{}}
\\\toprule
 \addlinespace
 \addlinespace
    Manifold                                                                                    & $f_0$ &  Range for $f_1$      &  Range for $g_2$    & Count \\ \midrule
 \addlinespace
 \addlinespace
 ${\mathbb R}{\bf P}^{\,3}$                                                                     &    11 &  $51\leq f_1\leq 52$  & $17\leq g_2\leq 18$ &     2 \\
 \addlinespace
\midrule
 \addlinespace
 \addlinespace
 $S^2\hbox{$\times\hspace{-1.62ex}\_\hspace{-.4ex}\_\hspace{.7ex}$}S^1$                         &    12 &  $60\leq f_1\leq 60$ & $22\leq g_2\leq 22$ &     2 \\
 \addlinespace
 ${\mathbb R}{\bf P}^{\,3}$                                                                     &    12 &  $60\leq f_1\leq 60$ & $22\leq g_2\leq 22$ &     4 \\
 \addlinespace
 $L(3,1)$                                                                                       &    12 &  $66\leq f_1\leq 66$ & $28\leq g_2\leq 28$ &     1 \\
 \addlinespace
 total                                                                                          &    12 &                      &                 &     7 \\
 \addlinespace
\midrule
 \addlinespace
 \addlinespace
 $S^2\hbox{$\times\hspace{-1.62ex}\_\hspace{-.4ex}\_\hspace{.7ex}$}S^1$                         &    13 &  $70\leq f_1\leq 75$ & $28\leq g_2\leq 33$ &     8 \\
 \addlinespace
 ${\mathbb R}{\bf P}^{\,3}$                                                                     &    13 &  $63\leq f_1\leq 63$ & $21\leq g_2\leq 21$ &     1 \\
 \addlinespace
 $L(3,1)$                                                                                       &    13 &  $70\leq f_1\leq 74$ & $28\leq g_2\leq 32$ &    72 \\
 \addlinespace
 total                                                                                          &    13 &                      &                 &    81 \\
 \addlinespace
\midrule
 \addlinespace
 \addlinespace
 $S^2\hbox{$\times\hspace{-1.62ex}\_\hspace{-.4ex}\_\hspace{.7ex}$}S^1$                         &    14 &  $75\leq f_1\leq 91$ & $29\leq g_2\leq 45$ &  6860 \\
 \addlinespace
 ${\mathbb R}{\bf P}^{\,3}$                                                                     &    14 &  $69\leq f_1\leq 70$ & $23\leq g_2\leq 24$ &     2 \\
 \addlinespace
 $L(3,1)$                                                                                       &    14 &  $75\leq f_1\leq 84$ & $29\leq g_2\leq 38$ &  7092 \\
 \addlinespace
 ${\mathbb R}{\bf P}^{\,2}\!\times S^1$                                                         &    14 &  $84\leq f_1\leq 91$ & $38\leq g_2\leq 45$ &  1011 \\
 \addlinespace
 $L(4,1)$                                                                                       &    14 &  $84\leq f_1\leq 90$ & $38\leq g_2\leq 44$ &   738 \\
 \addlinespace
 $L(5,2)$                                                                                       &    14 &  $86\leq f_1\leq 91$ & $40\leq g_2\leq 45$ &   121 \\
 \addlinespace
 total                                                                                          &    14 &                      &                 & 15824 \\
 \addlinespace
\midrule
 \addlinespace
 \addlinespace
 $S^3$                                                                                          &    15 &  $85\leq f_1\leq  94$ & $35\leq g_2\leq 44$ &      28 \\
  \addlinespace
 $S^2\hbox{$\times\hspace{-1.62ex}\_\hspace{-.4ex}\_\hspace{.7ex}$}S^1$                         &    15 &  $79\leq f_1\leq 104$ & $29\leq g_2\leq 54$ & 1500836 \\
 \addlinespace
 $S^2\!\times\!S^1$                                                                             &    15 &  $81\leq f_1\leq  97$ & $31\leq g_2\leq 47$ &      73 \\
 \addlinespace
 ${\mathbb R}{\bf P}^{\,3}$                                                                     &    15 &  $72\leq f_1\leq  94$ & $22\leq g_2\leq 44$ &      13 \\
 \addlinespace
 $L(3,1)$                                                                                       &    15 &  $81\leq f_1\leq  97$ & $31\leq g_2\leq 47$ &  240587 \\
 \addlinespace
 $(S^2\!\times\!S^1)^{\# 2}$                                                                    &    15 &  $91\leq f_1\leq 101$ & $41\leq g_2\leq 51$ &      84 \\
 \addlinespace
 $(S^2\hbox{$\times\hspace{-1.62ex}\_\hspace{-.4ex}\_\hspace{.7ex}$}S^1)^{\# 2}$                &    15 &  $95\leq f_1\leq 101$ & $45\leq g_2\leq 51$ &     144 \\
 \addlinespace
 ${\mathbb R}{\bf P}^{\,2}\!\times S^1$                                                         &    15 &  $88\leq f_1\leq 105$ & $38\leq g_2\leq 55$ & 3798307 \\
 \addlinespace
 $L(4,1)$                                                                                       &    15 &  $89\leq f_1\leq 101$ & $39\leq g_2\leq 51$ & 1968160 \\
 \addlinespace
 $L(5,2)$                                                                                       &    15 &  $90\leq f_1\leq 101$ & $40\leq g_2\leq 51$ &  504785 \\
 \addlinespace
 $(S^2\!\times\!S^1)\#{\mathbb R}{\bf P}^3$                                                     &    15 &  $91\leq f_1\leq  97$ & $41\leq g_2\leq 47$ &     238 \\
 \addlinespace
 $(S^2\hbox{$\times\hspace{-1.62ex}\_\hspace{-.4ex}\_\hspace{.7ex}$}S^1)\#{\mathbb R}{\bf P}^3$ &    15 &  $90\leq f_1\leq  99$ & $40\leq g_2\leq 49$ &    1913 \\
 \addlinespace
 $T^3$                                                                                          &    15 & $105\leq f_1\leq 105$ & $55\leq g_2\leq 55$ &       1 \\
 \addlinespace
 ${\mathbb R}{\bf P}^{\,3}\#\,{\mathbb R}{\bf P}^{\,3}$                                         &    15 &  $86\leq f_1\leq 102$ & $36\leq g_2\leq 52$ &  570885 \\
 \addlinespace
 $L(5,1)$                                                                                       &    15 &  $97\leq f_1\leq 102$ & $47\leq g_2\leq 52$ &    1314 \\
 \addlinespace
 $P_2=S^3/Q$                                                                                    &    15 &  $90\leq f_1\leq 102$ & $40\leq g_2\leq 52$ &   64475 \\
 \addlinespace
 $P_3$                                                                                          &    15 &  $97\leq f_1\leq 105$ & $47\leq g_2\leq 55$ &    1612 \\
 \addlinespace
 $P_4$                                                                                          &    15 & $104\leq f_1\leq 105$ & $54\leq g_2\leq 55$ &      20 \\
 \addlinespace
 $S^3/T^*$                                                                                      &    15 & $102\leq f_1\leq 102$ & $52\leq g_2\leq 52$ &       5 \\
 \addlinespace
 total                                                                                          &    15 &                       &                     & 8653480 \\
 \addlinespace
 \addlinespace
\bottomrule
\end{tabular}
\end{table}

To find all the (candidates for) $g_2$-irreducible triangulations with $f_0 \geq 16$ and $g_2 \leq 20$ 
a lower bound for $f_1$ was maintained during the construction of the 3-manifolds. 
When this lower bound became too large the search backtracked.  
The lower bound for $f_1$ was computed from the degrees of the finished vertices and lower bounds 
for the degrees of the other vertices. Theorem \ref{thm:lemma107} provides lower bounds for vertices 
which are neighbors of certain finished vertices.
Examples have been found which show that the lower bounds in Theorem \ref{thm:lemma107} are the best 
that can be obtained using just the necessary conditions of Theorems \ref{thm:lemma101}, \ref{thm:lemma102}, \ref{thm:lemma104}, and \ref{thm:lemma111}.

\begin{thm}  {\rm (Walkup \cite[10.7]{Walkup1970})} \label{thm:lemma107} 
Let $K$ be a $g_2$-irreducible triangulation and $(u,v)$ be an edge of $K$.  
\begin{itemize}
\addtolength{\itemsep}{-2.25mm}
\item If $L_v u$ is of type 9b(4'), 9b(4''), 9c(4) or 9c(4), then ${\deg}(v)\geq 7$.
\item If $L_v u$ is of type 8a(4), 8b(4), 9a(4), 9b(4) or 9d(4), then ${\deg}(v)\geq 8$.
\item If $L_v u$ is of type 7(4) or 9d(5), then ${\deg}(v)\geq 9$.
\item If $L_v u$ is of type 6(4), 8b(5), 9b(5), 9c(5) or 9c(5'), then ${\deg}(v)\geq 10$.
\item If $L_v u$ is of type 9b(6) or 9c(6) then ${\deg}(v)\geq 11$.
\item If $L_v u$ is of type 7(5), then ${\deg}(v)\geq 12$.
\item If $L_v u$ is of type 8a(6), then ${\deg}(v)\geq 14$.
\item If $L_v u$ is of type 9a(7), then ${\deg}(v)\geq 16$.
\item If the degree of $u$ is $10$ or $11$ and the degree of $(u,v)$ is $5$, then ${\deg}(v)\geq 8$.
\item If the degree of $u$ is $10$ and the degree of $(u,v)$ is $6$, then ${\deg}(v)\geq 10$.
\item If the degree of $u$ is $11$, $12$ or $13$ and the degree of $(u,v)$ is $6$, then ${\deg}(v)\geq 9$.
\item If the degree of $u$ is $11$, $12$, $13$, $14$ or $15$ and the degree of $(u,v)$ is $7$, then the degree of $v$ is at least $10$.
\item If the degree of $u$ is $11$ and the degree of $(u,v)$ is $6$, then ${\deg}(v)\geq 9$.
\item If the degree of $u$ is $11$ and the degree of $(u,v)$ is $7$, then ${\deg}(v)\geq 10$.
\item If the degree of $(u,v)$ is $d$, then the degree of $v$ is at least $d+2.$
\addtolength{\itemsep}{2.25mm}
\end{itemize}
\end{thm}

\textbf{Proof:}  
The last seven statements follow from the first eight.  
The first eight statements are from \cite[10.7]{Walkup1970}, \cite[11.2]{Walkup1970}, or \cite[11.4]{Walkup1970} or can be proved using the technique in the proof of \cite[10.7]{Walkup1970}.
Two of these results differ from \cite[10.7]{Walkup1970}.

To show that if $L_v u$ is of type 8b(5) then the degree of $v$  is at least 10 assume the degree of $v$ is less than 10.  
By \cite[11.1]{Walkup1970} $w$ (the bottom interior vertex in Figure~\ref{fig:WalkupFig3}-8b(5)) is in $W(u,v)$.  
$w$ is adjacent to four boundary vertices of 8b(5).  
For every type with the degree of $v$ less than 10 and five boundary vertices every interior vertex is adjacent to at least two boundary vertices.  
This contradicts \cite[10.6]{Walkup1970}.

If $L_v u$ is of type 8a(4) then the degree of $v$ may also be 8. 
Let $L_u v$ also be of type 8a(4),
let $W(u,v)$ be just the center vertex of Figure~\ref{fig:WalkupFig3}-8a(4), 
and identify the boundaries of $L_v u$ and $L_u v$ after rotating 
one copy of Figure~\ref{fig:WalkupFig3}-8a(4) a quarter of a turn.
\hfill$\Box$ 

\bigskip

The runs for $(f_0,f_1)=(16,73)$, $(16,74)$, $(17,78)$, and $(18,82)$ to search for 
potential $g_2$-irreducible triangulations produced no examples.
The run times were $1$, $4$, $64$, and $1000$ cpu-days, respectively.

\begin{thm}\label{thm:g2irred20}
The unique $g_2$-irreducible triangulation of ${\mathbb R}{\bf P}^3$ with $f=(11,51,80,40)$
and $g_2=17$ is the only $g_2$-irreducible triangulation of a $3$-manifold with $g_2 \leq 20$.
\end{thm}

\textbf{Proof:}
According to our enumeration, there are only two candidates for $g_2$-irreducible triangulations
with $g_2 \leq 20$; see Table~\ref{tbl:g2minimal_examples}.
Both candidates are triangulations of ${\mathbb R}{\bf P}^3$ with $11$ vertices.
One of the triangulations has $f$-vector $f=(11,51,80,40)$ and is the unique
$g_2$-irreducible triangulation of ${\mathbb R}{\bf P}^3$ by Theorem~\ref{thm:Walkup}.
The other candidate triangulation has $f$-vector $f=(11,52,82,41)$
and is therefore not $g_2$-minimal.\hfill$\Box$

\section{Examples of Triangulations and Upper Bounds 
      for \mathversion{bold}$\Gamma$ and $\Gamma^*$\mathversion{normal}}
\label{sec:examples}

Since $\Gamma(M)=\mbox{\rm min}\{\,g_2(K)\!\mid\! K\,\mbox{\rm\ is a triangulation of $M$\,}\}$ 
for any given $3$-mani\-fold~$M$, we get an upper bound $\Gamma(M)\leq g_2(K)$ for each
triangulation $K$ of $M$. Therefore, we are interested in ``small'' triangulations 
of the given manifold $M$. A standard procedure (cf.\ \cite{BjoernerLutz2000,Lutz2003bpre})
to obtain such small triangulations is to first construct any triangulation of $M$ of ``reasonable size'' 
and then to apply bistellar flips until a small or perhaps even vertex-minimal triangulation is reached.
The Tables~\ref{tbl:small_lens_spaces}--\ref{tbl:connected_sums_manifolds} list the $f$-vectors
of obtained triangulations along with resulting upper bounds on the respective $\Gamma$'s and $\Gamma^*$'s.
The triangulations that we found are available online at~\cite{Lutz_PAGE}.

According to Perelman's proof \cite{Perelman2003bpre}
of Thurston's Geometrization Conjecture \cite{Thurston1982},
every compact $3$-manifold can be decomposed canonically into geometric pieces,
which are modeled on one of eight model geometries. Six of the geometries, 
$S^3$ (spherical), $S^2\times{\mathbb R}^1$, $E^3$ (Euclidean), ${\rm Nil}$,
$H^2\times{\mathbb R}^1$, and $\widetilde{SL}(2,{\mathbb R})$,
yield Seifert manifolds, the other two geometries are ${\rm Sol}$ and $H^3$ (hyperbolic);
see~\cite{Scott1983} for a detailed discussion.

There are exactly four $3$-manifolds of geometry $S^2\times{\mathbb R}^1$
and ten flat  $3$-manifolds of geometry $E^3$, all other six geometries
give each rise to infinitely many $3$-mani\-folds.

The topological types of the $3$-manifolds modeled on the Seifert geometries 
are completely classified up to homeomorphism (cf.\ \cite{Orlik1972,Seifert1933}).
Moreover, it is possible to systematically construct triangulations of 
all Seifert manifolds; see \cite{BrehmLutz2002pre,Lutz2003bpre},
as well as \cite{Lutz_SEIFERT} for an implementation. For hyperbolic $3$-manifolds
it is unclear whether a complete classification can be obtained.
In 1982 Thurston \cite{Thurston1982} proved that almost every prime $3$-manifold 
is hyperbolic. Hyperbolic $3$-manifolds can be ordered with respect
to their hyperbolic volume. A census of 11,031 hyperbolic $3$-manifolds,
triangulated as pseudo-simplicial complexes with up to $30$ tetrahedra,
was obtained by Hodgson and Weeks \cite{HodgsonWeeks1994} by enumeration. 
For a census of all pseudo-simplicial triangulations of orientable 
and non-orientable $3$-manifolds with up to $10$ tetrahedra as well as 
for further references on pseudo-triangulation results, see Burton~\cite{Burton2007a}. 

In the Tables~\ref{tbl:small_lens_spaces}--\ref{tbl:small_nil_manifolds} we list manifolds 
of the six Seifert geometries, with the lens spaces of Table~\ref{tbl:small_lens_spaces}, 
the prism manifolds of Table~\ref{tbl:small_prism_manifolds},
and the three examples of Table~\ref{tbl:small_spherical_octa_dode}
of spherical geometry. 

The lens spaces $L(p,q)$, the prism spaces $P(r)$, the $Nil$ manifolds $\{ Oo,\hspace{1.1pt}1\mid b\}$,
and the products $M^2_{(+,g)}\!\times\!S^1$ and $M^2_{(-,g)}\!\times\!S^1$  have homology groups
\begin{eqnarray*}
H_*(L(p,q)) & = &({\mathbb Z},{\mathbb Z}_p,0,{\mathbb Z}), \\
H_*(P(r))   & = & \left\{
                   \begin{array}{l@{\hspace{5mm}}l}
                   ({\mathbb Z},{\mathbb Z}_2^{\,2},0,{\mathbb Z}),&\mbox{\rm $r$ even,}\\[1mm]
                   ({\mathbb Z},{\mathbb Z}_4,0,{\mathbb Z}),      &\mbox{\rm $r$ odd,}
                   \end{array}
                   \right.\\
H_*(\{ Oo,\hspace{1.1pt}1\mid b\}) & = & ({\mathbb Z},{\mathbb Z}^2\oplus {\mathbb Z}_b,{\mathbb Z}^2,{\mathbb Z}), \\
H_*(M^2_{(+,g)}\!\times\!S^1)      & = & ({\mathbb Z},{\mathbb Z}^{\,2g+1},{\mathbb Z}^{\,2g+1},{\mathbb Z}), \\
H_*(M^2_{(-,g)}\!\times\!S^1)      & = & ({\mathbb Z},{\mathbb Z}^g\oplus {\mathbb Z}_2,{\mathbb Z}^{\,g-1}\oplus {\mathbb Z}_2,0),
\end{eqnarray*}
respectively. Homology groups for the other examples are listed in the respective Tables.

Starting triangulations for the listed Seifert manifolds were either obtained 
by direct construction, as described in \cite{Lutz2003bpre}, or were produced 
with the program SEIFERT \cite{Lutz_SEIFERT} (cf.\ also \cite{BrehmLutz2002pre}). 
Small triangulations of these manifolds were already listed in \cite{Lutz2003bpre}. 
For a substantial number of the examples from \cite{Lutz2003bpre}
we were able to find yet smaller triangulations due to refinements 
of the bistellar flip technique and an increase of the number of
``rounds'' for the search.

The refined simulated annealing process consisted of three stages.
In the \emph{heating stage} we started with the best known triangulation of the
$3$-manifold of interest. The number of vertices was increased by half the number 
of vertices in the starting triangulation using only random $0$-moves;
i.e., moves were randomly chosen with $0$-moves, $1$-moves, $2$-moves, and $3$-moves
weighted by $1$, $0$, $0$, and $0$, respectively.
In the \emph{mixing stage} the heated triangulation was randomized without changing the number
of vertices; $10,000$ random $i$-moves were made with the four types of moves
weighted by $0$, $1$, $5$, and $0$.
In the subsequent \emph{cooling stage} the number of vertices was decreased whenever possible and
the number of edges was kept low; $100,000,000$ $i$-moves were made with the types
of moves weighted by $0$, $1$, $250$, and $\infty$.
The sequence of the mixing stage and the cooling stage was repeated ten times.
Any triangulation which had a smaller $f$-vector than the starting
triangulation was recorded.

\begin{table}
\small\centering
\defaultaddspace=0.2em
\caption{Lens spaces $L(p,q)$}\label{tbl:small_lens_spaces}
\begin{tabular}{@{}l@{\hspace{12mm}}l@{\hspace{12mm}}r@{\hspace{5mm}}r@{}}
\\\toprule
 \addlinespace
    Manifold         & Smallest known     & \multicolumn{2}{@{}l@{}}{Upper Bound for} \\
 \addlinespace
                     & $f$-vector         & \mbox{}\hspace{7mm}$\Gamma^*$ & $\Gamma$  \\ \midrule
 \addlinespace
 \addlinespace
 $L(1,1)=S^3$        & (5,10,10,5)        &  0 &  0 \\
 \addlinespace
 $L(2,1)={\mathbb R}{\bf P}^3$   & (11,51,80,40) & 17 & 17 \\
 \addlinespace
 $L(3,1)$            & (12,66,108,54)     & 28 & 28 \\
 \addlinespace
 $L(4,1)$            & (14,84,140,70)     & 38 & 38 \\
 \addlinespace
 $L(5,1)$            & (15,97,164,82)     & 47 & 47 \\
 \addlinespace
 $L(6,1)$            & (16,110,188,94)    & 56 & 56 \\
 \addlinespace
 $L(7,1)$            & (17,123,212,106)   & 67 & 65 \\
 \addlinespace
 $L(8,1)$            & (17,130,226,113)   & 72 & 72 \\
 \addlinespace
 $L(9,1)$            & (18,143,252,126)   & 81 & 81 \\
 \addlinespace
 $L(10,1)$           & (19,155,272,136)   & 92 & 89 \\
 \addlinespace
\midrule
 \addlinespace
 \addlinespace
 $L(5,2)$            & (14,86,144,72)     & 40 & 40 \\
 \addlinespace
 $L(7,2)$            & (16,104,176,88)    & 56 & 50 \\
 \addlinespace
 $L(8,3)$            & (16,106,180,90)    & 56 &    \\
                     & [(17,109,184,92)]  &    & 51 \\
 \addlinespace
 $L(9,2)$            & (16,114,196,98)    & 60 &    \\    
                     & [(17,116,198,99)]  &    & 58 \\
 \addlinespace
 $L(10,3)$           & (17,118,202,101)   & 67 &    \\   
                     & [(18,121,206,103)] &    & 59 \\
 \addlinespace
 \addlinespace
\bottomrule
\end{tabular}
\end{table}

\begin{table}
\small\centering
\defaultaddspace=0.2em
\caption{Prism manifolds}\label{tbl:small_prism_manifolds}
\begin{tabular}{@{}l@{\hspace{12mm}}l@{\hspace{12mm}}r@{\hspace{5mm}}r@{}}
\\\toprule
 \addlinespace
    Manifold         & Smallest known     & \multicolumn{2}{@{}l@{}}{Upper Bound for} \\
 \addlinespace
                     & $f$-vector         & \mbox{}\hspace{7mm}$\Gamma^*$ & $\Gamma$  \\ \midrule
 \addlinespace
 \addlinespace
 $P_2=S^3/Q$, cube space  & (15,90,150,75)     &  46 &  40 \\
 \addlinespace
 $P_3$                    & (15,97,164,82)     &  47 &  47 \\
 \addlinespace
 $P_4$                    & (15,104,178,89)    &  54 &  54 \\
 \addlinespace
 $P_5$                    & (17,122,210,105)   &  67 &  64 \\
 \addlinespace
 $P_6$                    & (17,130,226,113)   &  72 &  72 \\
 \addlinespace
 $P_7$                    & (18,143,250,125)   &  81 &  81 \\
 \addlinespace
 $P_8$                    & (19,155,272,136)   &  92 &  89 \\
 \addlinespace
 $P_9$                    & (19,163,288,144)   &  97 &  97 \\
 \addlinespace
 $P_{10}$                 & (20,175,310,155)   & 106 & 105 \\
 \addlinespace
 \addlinespace
\bottomrule
\end{tabular}
\end{table}

\begin{table}
\small\centering
\defaultaddspace=0.15em
\caption{The spherical octahedral, truncated cube, and dodecahedral space}\label{tbl:small_spherical_octa_dode}
\begin{tabular*}{\linewidth}{@{\extracolsep{\fill}}lll@{\hspace{-4mm}}rr@{}}
\\\toprule
 \addlinespace
    Manifold         & Homology & Smallest known     & \multicolumn{2}{@{}l@{}}{\hfill Upper Bound for} \\
 \addlinespace
                     &          & $f$-vector         & \mbox{}\hspace{15mm}$\Gamma^*$ & $\Gamma$  \\ \midrule
 \addlinespace
 \addlinespace
  $S^3/T^*$          & $({\mathbb Z},{\mathbb Z}_3,0,{\mathbb Z})$ & (15,102,174,87)   & 52 &    \\
                     &                                             & [(16,104,176,88)] &    & 50 \\
 \addlinespace
  $S^3/O^*$          & $({\mathbb Z},{\mathbb Z}_2,0,{\mathbb Z})$ & (16,109,186,93)   & 56 & 55 \\
 \addlinespace
 $S^3/I^*=\Sigma(2,3,5)$, \\
 Poincar\'e $3$-sphere & $({\mathbb Z},0,0,{\mathbb Z})$           & (16,106,180,90)   & 56 & 52 \\
 \addlinespace
 \addlinespace
\bottomrule
\end{tabular*}
\end{table}

\begin{table}
\small\centering
\defaultaddspace=0.15em
\caption{($S^2\times{\mathbb R}$)-spaces}\label{tbl:S2xR_spaces}
\begin{tabular*}{\linewidth}{@{\extracolsep{\fill}}lll@{\hspace{-4mm}}rr@{}}
\\\toprule
 \addlinespace
    Manifold         & Homology & Smallest known     & \multicolumn{2}{@{}l@{}}{\hfill Upper Bound for} \\
 \addlinespace
                     &          & $f$-vector         & \mbox{}\hspace{17mm}$\Gamma^*$ & $\Gamma$  \\ \midrule
 \addlinespace
 \addlinespace
 $S^2\hbox{$\times\hspace{-1.62ex}\_\hspace{-.4ex}\_\hspace{.7ex}$}S^1$ & $({\mathbb Z},{\mathbb Z},{\mathbb Z}_2,0)$ & (9,36,54,27)   & 10 & 10 \\
 \addlinespace
 $S^2\!\times\!S^1$                                      & $({\mathbb Z},{\mathbb Z},{\mathbb Z},{\mathbb Z})$        & (10,40,60,30)  & 11 & 10 \\
 \addlinespace
 ${\mathbb R}{\bf P}^{\,2}\!\times S^1$                  & $({\mathbb Z},{\mathbb Z}\oplus{\mathbb Z}_2,{\mathbb Z}_2,0)$ & (14,84,140,70)  & 38 & 38 \\
 \addlinespace
 ${\mathbb R}{\bf P}^{\,3}\#\,{\mathbb R}{\bf P}^{\,3}$  & $({\mathbb Z},{\mathbb Z}_2^{\,2},0,{\mathbb Z})$ & (15,86,142,71)   & 46 &    \\
                                                         &                                                   & [(16,89,146,73)] &    &    \\
                                                         &                                                   & [(18,96,156,78)] &    & 34 \\
 \addlinespace
 \addlinespace
\bottomrule
\end{tabular*}
\end{table}

\begin{table}
\small\centering
\defaultaddspace=0.15em
\caption{Flat manifolds}\label{tbl:small_flat_manifolds}
\begin{tabular*}{\linewidth}{@{\extracolsep{\fill}}lll@{\hspace{-3mm}}rr@{}}
\\\toprule
 \addlinespace
    Manifold      & Homology & Smallest known     & \multicolumn{2}{@{}l@{}}{\hfill Upper Bound for} \\
 \addlinespace
                  &          & $f$-vector         & \mbox{}\hspace{14mm}$\Gamma^*$ & $\Gamma$  \\ \midrule
 \addlinespace
 \addlinespace
 $T^3$            & $({\mathbb Z},{\mathbb Z}^3,{\mathbb Z}^3,{\mathbb Z})$                           & (15,105,180,90)    & 55 &    \\
                  &                                                                                   & [(16,108,184,92)]  &    & 54 \\
 \addlinespace
 $G_2$            & $({\mathbb Z},{\mathbb Z}\oplus{\mathbb Z}_2^{\,2},{\mathbb Z},{\mathbb Z})$      & (16,116,200,100)   & 62 &    \\
                  &                                                                                   & [(17,118,202,101)] &    & 60 \\
 \addlinespace
 $G_3$            & $({\mathbb Z},{\mathbb Z}\oplus{\mathbb Z}_3,{\mathbb Z},{\mathbb Z})$            & (17,117,200,100)   & 67 & 59 \\
 \addlinespace
 $G_4$            & $({\mathbb Z},{\mathbb Z}\oplus{\mathbb Z}_2,{\mathbb Z},{\mathbb Z})$            & (16,115,198,99)    & 61 & 61 \\
 \addlinespace
 $G_5$            & $({\mathbb Z},{\mathbb Z},{\mathbb Z},{\mathbb Z})$                               & (16,112,192,96)    & 58 &    \\ 
                  &                                                                                   & [(17,115,196,98)]  &    & 57 \\
 \addlinespace
 $G_6$            & $({\mathbb Z},{\mathbb Z}_4^{\,2},0,{\mathbb Z})$                                 & (17,124,214,107)   & 67 & 66 \\
 \addlinespace
\midrule
 \addlinespace
 \addlinespace
 $K\times S^1$    & $({\mathbb Z},{\mathbb Z}^2\oplus{\mathbb Z}_2,{\mathbb Z}\oplus{\mathbb Z}_2,0)$ & (16,115,198,99)    & 61 &    \\
                  &                                                                                   & [(17,118,202,101)] &    & 60 \\
 \addlinespace
 $B_2$            & $({\mathbb Z},{\mathbb Z}^2,{\mathbb Z}\oplus{\mathbb Z}_2,0)$                    & (16,110,188,94)    & 56 & 56 \\
 \addlinespace
 $B_3$            & $({\mathbb Z},{\mathbb Z}\oplus{\mathbb Z}_2^{\,2},{\mathbb Z}_2,0)$              & (17,119,204,102)   & 67 & 61 \\
 \addlinespace
 $B_4$            & $({\mathbb Z},{\mathbb Z}\oplus{\mathbb Z}_4,{\mathbb Z}_2,0)$                    & (17,117,200,100)   & 67 & 59 \\
 \addlinespace
 \addlinespace
\bottomrule
\end{tabular*}
\end{table}

\begin{table}
\small\centering
\defaultaddspace=0.2em
\caption{($H^2\times{\mathbb R}$)-spaces}\label{tbl:H2xR_spaces}
\begin{tabular}{@{}l@{\hspace{12mm}}l@{\hspace{12mm}}r@{\hspace{5mm}}r@{}}
\\\toprule
 \addlinespace
    Manifold         & Smallest known     & \multicolumn{2}{@{}l@{}}{Upper Bound for} \\
 \addlinespace
                     & $f$-vector         & \mbox{}\hspace{7mm}$\Gamma^*$ & $\Gamma$  \\ \midrule
 \addlinespace
 \addlinespace
 $M^2_{(+,2)}\!\times\!S^1$    & (20,168,296,148)   & 106 &  98 \\
 \addlinespace
 $M^2_{(+,3)}\!\times\!S^1$    & (22,210,376,188)   & 137 & 132 \\
 \addlinespace
 $M^2_{(+,4)}\!\times\!S^1$    & (24,256,464,232)   & 172 & 170 \\
 \addlinespace
 $M^2_{(+,5)}\!\times\!S^1$    & (26,299,546,273)   & 211 & 205 \\
 \addlinespace
\midrule
 \addlinespace
 \addlinespace
 $M^2_{(-,3)}\!\times\!S^1$    & (18,141,246,123)   &  79 &  79 \\
 \addlinespace
 $M^2_{(-,4)}\!\times\!S^1$    & (19,163,288,144)   &  97 &     \\
                               & [(20,166,292,146)] &     &  96 \\
 \addlinespace
 $M^2_{(-,5)}\!\times\!S^1$    & (21,190,338,169)   & 121 & 116 \\
 \addlinespace
 $M^2_{(-,6)}\!\times\!S^1$    & (22,212,380,190)   & 137 & 134 \\
 \addlinespace
 $M^2_{(-,7)}\!\times\!S^1$    & (23,234,422,211)   & 154 & 152 \\
 \addlinespace
 $M^2_{(-,8)}\!\times\!S^1$    & (24,256,464,232)   & 172 &     \\
                               & [(25,259,468,234)] &     & 169 \\
 \addlinespace
 $M^2_{(-,9)}\!\times\!S^1$    & (25,277,504,252)   & 191 & 187 \\
 \addlinespace
 $M^2_{(-,10)}\!\times\!S^1$   & (26,296,540,270)   & 211 & 202 \\
 \addlinespace
 \addlinespace
\bottomrule
\end{tabular}
\end{table}

\begin{table}
\small\centering
\defaultaddspace=0.2em
\caption{Seifert homology spheres of geometry $SL(2,{\mathbb Z})$}\label{tbl:small_SL2Z_spheres}
\begin{tabular}{@{}l@{\hspace{12mm}}l@{\hspace{12mm}}r@{\hspace{5mm}}r@{}}
\\\toprule
 \addlinespace
    Manifold         & Smallest known     & \multicolumn{2}{@{}l@{}}{Upper Bound for} \\
 \addlinespace
                     & $f$-vector         & \mbox{}\hspace{7mm}$\Gamma^*$ & $\Gamma$  \\ \midrule
 \addlinespace
 \addlinespace
 $\Sigma(2,3,7)$         & (16,117,202,101)   &  63 &  63 \\
 \addlinespace
 $\Sigma(2,5,7)$         & (18,138,240,120)   &  79 &  76 \\
 \addlinespace
 $\Sigma(3,4,5)$         & (18,139,242,121)   &  79 &  77 \\
 \addlinespace
 $\Sigma(3,4,7)$         & (18,151,266,133)   &  89 &     \\
                         & [(19,153,268,134)] &     &     \\
                         & [(20,156,272,136)] &     &  86 \\
 \addlinespace
 $\Sigma(3,5,7)$         & (20,171,302,151)   & 106 & 101 \\
 \addlinespace
 $\Sigma(4,5,7)$         & (20,177,314,157)   & 107 &     \\
                         & [(21,179,316,158)] &     & 105 \\
 \addlinespace
 \addlinespace
\bottomrule
\end{tabular}
\end{table}

\begin{table}
\small\centering
\defaultaddspace=0.2em
\caption{$Nil$ manifolds}\label{tbl:small_nil_manifolds}
\begin{tabular}{@{}l@{\hspace{12mm}}l@{\hspace{12mm}}r@{\hspace{5mm}}r@{}}
\\\toprule
 \addlinespace
    Manifold         & Smallest known     & \multicolumn{2}{@{}l@{}}{Upper Bound for} \\
 \addlinespace
                     & $f$-vector         & \mbox{}\hspace{7mm}$\Gamma^*$ & $\Gamma$  \\ \midrule
 \addlinespace
 \addlinespace
 $\{Oo,1\mid 1\}$        & (16,113,194,97)  & 59 & 59 \\
 \addlinespace
 $\{Oo,1\mid 2\}$        & (17,120,206,103) & 67 & 62 \\
 \addlinespace
 $\{Oo,1\mid 3\}$        & (17,125,216,108) & 67 & 67 \\
 \addlinespace
 $\{Oo,1\mid 4\}$        & (17,130,226,113) & 72 & 72 \\
 \addlinespace
 $\{Oo,1\mid 5\}$        & (18,142,248,124) & 80 & 80 \\
 \addlinespace
 \addlinespace
\bottomrule
\end{tabular}
\end{table}

The examples of $3$-manifolds of ${\rm Sol}$ geometry are either torus or Klein bottle
bundles over $S^1$ or are composed of two twisted $I$-bundles over the torus or the Klein bottle;
cf.\ Hempel and Jaco \cite{HempelJaco1972} and Scott \cite{Scott1983}. However, a topological 
classification of the individual examples seems not to be known, which kept us from 
providing explicit examples of this geometry.

The Hodgson--Weeks census from 1994 \cite{HodgsonWeeks1994} is still the main source for 
(pseudo-sim\-pli\-cial triangulations of) closed hyperbolic $3$-manifold. Most of the 
census examples are orientable. The orientable example \texttt{or\_0.94270736}
of smallest listed hyperbolic volume $0.94270736$ is called the Weeks manifold.
It has recently been proved by Gabai, Meyerhoff, and Milley~\cite{GabaiMeyerhoffMilley2007pre}
that the Weeks manifold with homology $({\mathbb Z},{\mathbb Z}_5^{\,2},0,{\mathbb Z})$ 
 has smallest possible volume among all orientable hyperbolic $3$-manifolds.

The hyperbolic census data is accessible via the \texttt{SnapPea}-package~\cite{SnapPea}
of Weeks, via the \texttt{Regina}-package \cite{regina} of Burton (cf.\ also \cite{Burton2004}), 
or online via \url{http://regina.sourceforge.net/data.html}. For our purposes
it was necessary to first turn the pseudo-simplicial
complexes, consisting of a set of tetrahedra with gluing information for the 
boundaries, into proper simplicial complexes. The listed pseudo-triangulations
all have only one vertex and between $9$--$30$ tetrahedra. If the number of
starting tetrahedra is $\mbox{\rm ntet}$, the second barycentric subdivision
is a proper simplicial complex with $24^2\cdot\mbox{\rm ntet}$ tetrahedra.
The desired small triangulations are then obtained via bistellar flips.
For example, the second barycentric subdivision of the Weeks manifold 
with $\mbox{\rm ntet}=9$ has $f$-vector $f=(940,6124,10368,5184)$.
In this case, the smallest triangulation of the Weeks manifold that 
we found has $f=(18,141,246,123)$; see Table~\ref{tbl:small_hyper}.

\begin{table}
\small\centering
\defaultaddspace=0.2em
\caption{The first twenty hyperbolic $3$-manifolds from the Hodgson--Weeks census
         and the Weber--Seifert hyperbolic dodecahedral space.}\label{tbl:small_hyper}
\begin{tabular*}{\linewidth}{@{\extracolsep{\fill}}lll@{\hspace{-1mm}}rr@{}}
\\\toprule
 \addlinespace
    Manifold         & Homology & Smallest known     & \multicolumn{2}{@{}l@{}}{\hfill Upper Bound for} \\
 \addlinespace
                     &          & $f$-vector         & \mbox{}\hspace{11mm}$\Gamma^*$ & $\Gamma$ \\ \midrule
 \addlinespace
 \addlinespace
  \texttt{or\_0.94270736}    & $({\mathbb Z},{\mathbb Z}_5^{\,2},0,{\mathbb Z})$ & (18,139,242,121)   & 79 &    \\
                             &                                                   & [(19,142,246,123)] &    & 76 \\
 \addlinespace
  \texttt{or\_0.98136883}    & $({\mathbb Z},{\mathbb Z}_5,0,{\mathbb Z})$       & (18,135,234,117)   & 79 & 73 \\
 \addlinespace
  \texttt{or\_1.01494161}    & $({\mathbb Z},{\mathbb Z}_3+{\mathbb Z}_6,0,{\mathbb Z})$ & (18,135,234,117) & 79 & 73 \\
 \addlinespace
  \texttt{or\_1.26370924}    & $({\mathbb Z},{\mathbb Z}_5^{\,2},0,{\mathbb Z})$ & (18,149,262,131)   & 87 &    \\
                             &                                                   & [(19,150,262,131)] &    & 84 \\
 \addlinespace
  \texttt{or\_1.28448530}    & $({\mathbb Z},{\mathbb Z}_6,0,{\mathbb Z})$       & (18,139,242,121)   & 79 & 77 \\
 \addlinespace
  \texttt{or\_1.39850888}    & $({\mathbb Z},0,0,{\mathbb Z})$                   & (18,140,244,122)   & 79 &    \\
                             &                                                   & [(19,143,248,124)] &    & 77 \\
 \addlinespace
  \texttt{or\_1.41406104\_a} & $({\mathbb Z},{\mathbb Z}_6,0,{\mathbb Z})$       & (18,136,236,118)   & 79 & 74 \\
 \addlinespace
  \texttt{or\_1.41406104\_b} & $({\mathbb Z},{\mathbb Z}_{10},0,{\mathbb Z})$    & (18,145,254,127)   & 83 &    \\
                             &                                                   & [(19,147,256,128)] &    & 81 \\ 
 \addlinespace
  \texttt{or\_1.42361190}    & $({\mathbb Z},{\mathbb Z}_{35},0,{\mathbb Z})$    & (19,153,268,134)   & 92 & 87 \\
 \addlinespace
  \texttt{or\_1.44069901}    & $({\mathbb Z},{\mathbb Z}_3,0,{\mathbb Z})$       & (18,141,246,123)   & 79 & 79 \\
 \addlinespace
  \texttt{or\_1.46377664}    & $({\mathbb Z},{\mathbb Z}_7,0,{\mathbb Z})$       & (18,148,260,130)   & 86 &    \\
                             &                                                   & [(19,150,262,131)] &    & 84 \\ 
 \addlinespace
  \texttt{or\_1.52947733}    & $({\mathbb Z},{\mathbb Z}_5,0,{\mathbb Z})$       & (18,144,252,126)   & 82 &    \\
                             &                                                   & [(19,146,254,127)] &    & 80 \\ 
 \addlinespace
  \texttt{or\_1.54356891\_a} & $({\mathbb Z},{\mathbb Z}_{35},0,{\mathbb Z})$    & (19,152,266,133)   & 92 &    \\
                             &                                                   & [(21,159,276,138)] &    & 85 \\ 
 \addlinespace
  \texttt{or\_1.54356891\_b} & $({\mathbb Z},{\mathbb Z}_{21},0,{\mathbb Z})$    & (18,144,252,126)   & 82 &    \\
                             &                                                   & [(19,147,256,128)] &    & 81 \\ 
 \addlinespace
  \texttt{or\_1.58316666\_a} & $({\mathbb Z},{\mathbb Z}_{21},0,{\mathbb Z})$    & (18,140,244,122)   & 79 & 78 \\
 \addlinespace
  \texttt{or\_1.58316666\_b} & $({\mathbb Z},{\mathbb Z}_3+{\mathbb Z}_9,0,{\mathbb Z})$ & (18,144,252,126) & 82 &    \\
                             &                                                   & [(19,147,256,128)] &    &    \\ 
                             &                                                   & [(20,150,260,130)] &    & 80 \\ 
 \addlinespace
  \texttt{or\_1.58864664\_a} & $({\mathbb Z},{\mathbb Z}_{30},0,{\mathbb Z})$    & (19,151,264,132)   & 92 & 85 \\
 \addlinespace
  \texttt{or\_1.58864664\_b} & $({\mathbb Z},{\mathbb Z}_{30},0,{\mathbb Z})$    & (19,159,280,140)   & 93 &     \\
                             &                                                   & [(20,162,284,142)] &    & 92 \\ 
 \addlinespace
  \texttt{or\_1.64960972}    & $({\mathbb Z},{\mathbb Z}_{15},0,{\mathbb Z})$    & (18,147,258,129)   & 85 &    \\
                             &                                                   & [(19,150,262,131)] &    & 84 \\ 
 \addlinespace
  \texttt{or\_1.75712603}    & $({\mathbb Z},{\mathbb Z}_7,0,{\mathbb Z})$       & (18,140,244,122)   & 79 & 78 \\ 
 \addlinespace
\midrule
 \addlinespace
 \addlinespace
  hyperb.\ dodec.\ space     & $({\mathbb Z},{\mathbb Z}_5^{\,3},0,{\mathbb Z})$ & (21,190,338,169) & 121 & 116 \\
 \addlinespace
 \addlinespace
\bottomrule
\end{tabular*}
\end{table}

A well-known example of a hyperbolic $3$-manifold is the Weber--Seifert hyperbolic dodecahedral space
with homology $({\mathbb Z},{\mathbb Z}_5^{\,3},0,{\mathbb Z})$. Our smallest triangulation
of this manifold has $f=(21,193,344,172)$, which is close to the $18$ vertices
of the Weeks manifold. Nevertheless, the Weber--Seifert hyperbolic dodecahedral space
does not appear in the Hodgson--Weeks census (there is no manifold with this homology
in the census). In fact, sixteen of the first twenty examples from the
census have triangulations as simplicial complexes with $18$ vertices, 
the remaining four with $19$ vertices; see Table~\ref{tbl:small_hyper}.
This seems to indicate that perhaps most of the 11,031 census examples have triangulations
as proper simplicial complexes with $18$--$21$ vertices. 

\begin{conj}
At least $18$ vertices are needed to triangulate a hyperbolic $3$-manifold
as a simplicial complex.
\end{conj}

It was proved by Brehm and Swiatkowski \cite{BrehmSwiatkowski1993} that
the number of non-ho\-meo\-mor\-phic lens spaces that can be triangulated
with $n$ vertices grows exponentially with $n$. 

In contrast to the many triangulations that we expect from $18$ vertices on,
the list of $3$-manifolds that can be triangulated with at most $17$ vertices
will be comparably short. In \cite{Lutz2005bpre}, $27$ different $3$-manifolds 
were described that can be triangulated with up to $15$ vertices. With our 
improved bistellar flip techniques we were able to find further $6$.

\begin{thm}
There are at least $33$ different $3$-manifolds that can be triangulated with up to $15$ vertices.
These examples are:
\begin{tabbing}
mm\= mmmmm \= \kill
\> $n=5$:  \> $S^3$, \\[.25mm]
\> $n=9$:  \> $S^2\hbox{$\times\hspace{-1.62ex}\_\hspace{-.4ex}\_\hspace{.7ex}$}S^1$, \\[.25mm]
\> $n=10$: \> $S^2\!\times\!S^1$, \\[.25mm]
\> $n=11$: \> ${\mathbb R}{\bf P}^3$, \\[.25mm]
\> $n=12$: \> $L(3,1)$, $(S^2\!\times\!S^1)^{\# 2}$,\, $(S^2\hbox{$\times\hspace{-1.62ex}\_\hspace{-.4ex}\_\hspace{.7ex}$}S^1)^{\# 2}$, \\[.25mm]
\> $n=13$: \> $(S^2\!\times\!S^1)^{\# 3}$,\, $(S^2\hbox{$\times\hspace{-1.62ex}\_\hspace{-.4ex}\_\hspace{.7ex}$}S^1)^{\# 3}$, \\[.25mm]
\> $n=14$: \> ${\mathbb R}{\bf P}^{\,2}\!\times S^1$,\, $L(4,1)$,\, $L(5,2)$,
              $(S^2\!\times\!S^1)^{\# 4}$,\, $(S^2\hbox{$\times\hspace{-1.62ex}\_\hspace{-.4ex}\_\hspace{.7ex}$}S^1)^{\# 4}$,\\
\>         \> $(S^2\!\times\!S^1)\#{\mathbb R}{\bf P}^3$, $(S^2\hbox{$\times\hspace{-1.62ex}\_\hspace{-.4ex}\_\hspace{.7ex}$}S^1)\#{\mathbb R}{\bf P}^3$,\\
\>         \> $(S^2\!\times\!S^1)^{\# 2}\#{\mathbb R}{\bf P}^3$, $(S^2\hbox{$\times\hspace{-1.62ex}\_\hspace{-.4ex}\_\hspace{.7ex}$}S^1)^{\# 2}\#{\mathbb R}{\bf P}^3$, \\[.25mm]
\> $n=15$: \> ${\mathbb R}{\bf P}^{\,3}\#\,{\mathbb R}{\bf P}^{\,3}$,\, $L(5,1)$,\, $S^3/Q$,\, $P_3$,\, $P_4$,\, $S^3/T^*$,\, $T^3$, \\
\>         \> $(S^2\!\times\!S^1)^{\# 5}$,\, $(S^2\hbox{$\times\hspace{-1.62ex}\_\hspace{-.4ex}\_\hspace{.7ex}$}S^1)^{\# 5}$,\, 
              $(S^2\!\times\!S^1)^{\# 3}\#{\mathbb R}{\bf P}^3$,\, $(S^2\hbox{$\times\hspace{-1.62ex}\_\hspace{-.4ex}\_\hspace{.7ex}$}S^1)^{\# 3}\#{\mathbb R}{\bf P}^3$,\\
\>         \> $(S^2\!\times\!S^1)\#L(3,1)$, $(S^2\hbox{$\times\hspace{-1.62ex}\_\hspace{-.4ex}\_\hspace{.7ex}$}S^1)\#L(3,1)$,\\
\>         \> $(S^2\!\times\!S^1)^{\# 2}\#L(3,1)$,\, $(S^2\hbox{$\times\hspace{-1.62ex}\_\hspace{-.4ex}\_\hspace{.7ex}$}S^1)^{\# 2}\#L(3,1)$.
\end{tabbing}
\end{thm}
It is conjectured in \cite{Lutz2005bpre} that this list is complete up to $13$ vertices. 
The particular examples are listed in 
Tables~\ref{tbl:small_lens_spaces}--\ref{tbl:small_nil_manifolds} (Seifert manifolds) and in 
Table~\ref{tbl:connected_sums_manifolds} (connected sums of Seifert manifolds).

\begin{table}
\small\centering
\defaultaddspace=0.2em
\caption{Connected sums}\label{tbl:connected_sums_manifolds}
\begin{tabular}{@{}l@{\hspace{8mm}}l@{\hspace{8mm}}r@{\hspace{5mm}}r@{}}
\\\toprule
 \addlinespace
    Manifold         & Smallest known     & \multicolumn{2}{@{}l@{}}{Upper Bound for} \\
 \addlinespace
                     & $f$-vector         & \mbox{}\hspace{7mm}$\Gamma^*$ & $\Gamma$  \\ \midrule
 \addlinespace
 \addlinespace
 $(S^2\!\times\!S^1)^{\# 0}=S^3$                                                                                   & (5,10,10,5)       &   0 &   0 \\
 \addlinespace
 $(S^2\hbox{$\times\hspace{-1.62ex}\_\hspace{-.4ex}\_\hspace{.7ex}$}S^1)^{\# 1}$                                   & (9,36,54,27)      &  10 &  10 \\
 \addlinespace
 $(S^2\!\times\!S^1)^{\# 1}$                                                                                       & (10,40,60,30)     &  11 &  10 \\
 \addlinespace
\midrule
 \addlinespace
 \addlinespace
  $(S^2\!\times\!S^1)^{\# 2}$,\, $(S^2\hbox{$\times\hspace{-1.62ex}\_\hspace{-.4ex}\_\hspace{.7ex}$}S^1)^{\# 2}$   & (12,58,92,46)     &  22 &  20 \\
 \addlinespace
  $(S^2\!\times\!S^1)^{\# 3}$,\, $(S^2\hbox{$\times\hspace{-1.62ex}\_\hspace{-.4ex}\_\hspace{.7ex}$}S^1)^{\# 3}$   & (13,72,118,59)    &  30 &  30 \\
 \addlinespace
  $(S^2\!\times\!S^1)^{\# 4}$,\, $(S^2\hbox{$\times\hspace{-1.62ex}\_\hspace{-.4ex}\_\hspace{.7ex}$}S^1)^{\# 4}$   & (14,86,144,72)    &  40 &  40 \\
 \addlinespace
  $(S^2\!\times\!S^1)^{\# 5}$,\, $(S^2\hbox{$\times\hspace{-1.62ex}\_\hspace{-.4ex}\_\hspace{.7ex}$}S^1)^{\# 5}$   & (15,100,170,85)   &  50 &  50 \\
 \addlinespace
  $(S^2\!\times\!S^1)^{\# 6}$,\, $(S^2\hbox{$\times\hspace{-1.62ex}\_\hspace{-.4ex}\_\hspace{.7ex}$}S^1)^{\# 6}$   & (16,114,196,98)   &  60 &  60 \\
 \addlinespace
  $(S^2\!\times\!S^1)^{\# 7}$,\, $(S^2\hbox{$\times\hspace{-1.62ex}\_\hspace{-.4ex}\_\hspace{.7ex}$}S^1)^{\# 7}$   & (17,128,222,111)  &  70 &  70 \\
 \addlinespace
  $(S^2\!\times\!S^1)^{\# 8}$,\, $(S^2\hbox{$\times\hspace{-1.62ex}\_\hspace{-.4ex}\_\hspace{.7ex}$}S^1)^{\# 8}$   & (18,142,248,124)  &  80 &  80 \\
 \addlinespace
  $(S^2\!\times\!S^1)^{\# 9}$,\, $(S^2\hbox{$\times\hspace{-1.62ex}\_\hspace{-.4ex}\_\hspace{.7ex}$}S^1)^{\# 9}$   & (19,156,274,137)  &  92 &  90 \\
 \addlinespace
  $(S^2\!\times\!S^1)^{\# 10}$,\, $(S^2\hbox{$\times\hspace{-1.62ex}\_\hspace{-.4ex}\_\hspace{.7ex}$}S^1)^{\# 10}$ & (19,166,294,147)  & 100 & 100 \\
 \addlinespace
  $(S^2\!\times\!S^1)^{\# 11}$,\, $(S^2\hbox{$\times\hspace{-1.62ex}\_\hspace{-.4ex}\_\hspace{.7ex}$}S^1)^{\# 11}$ & (20,180,320,160)  & 110 & 110 \\
 \addlinespace
  $(S^2\!\times\!S^1)^{\# 12}$,\, $(S^2\hbox{$\times\hspace{-1.62ex}\_\hspace{-.4ex}\_\hspace{.7ex}$}S^1)^{\# 12}$ & (21,194,346,173)  & 121 & 120 \\
 \addlinespace
  $(S^2\!\times\!S^1)^{\# 13}$,\, $(S^2\hbox{$\times\hspace{-1.62ex}\_\hspace{-.4ex}\_\hspace{.7ex}$}S^1)^{\# 13}$ & (22,208,372,186)  & 137 & 130 \\
 \addlinespace
  $(S^2\!\times\!S^1)^{\# 14}$,\, $(S^2\hbox{$\times\hspace{-1.62ex}\_\hspace{-.4ex}\_\hspace{.7ex}$}S^1)^{\# 14}$ & (22,218,392,196)  & 140 & 140 \\
 \addlinespace
  $(S^2\!\times\!S^1)^{\# 15}$,\, $(S^2\hbox{$\times\hspace{-1.62ex}\_\hspace{-.4ex}\_\hspace{.7ex}$}S^1)^{\# 15}$ & (23,232,418,209)  & 154 & 150 \\
 \addlinespace
\midrule
 \addlinespace
 \addlinespace
 $(S^2\!\times\!S^1)\#{\mathbb R}{\bf P}^3$,\, $(S^2\hbox{$\times\hspace{-1.62ex}\_\hspace{-.4ex}\_\hspace{.7ex}$}S^1)\#{\mathbb R}{\bf P}^3$
         & (14,73,118,59) & 37 & 27 \\
 \addlinespace
 $(S^2\!\times\!S^1)^{\# 2}\#{\mathbb R}{\bf P}^3$,\, $(S^2\hbox{$\times\hspace{-1.62ex}\_\hspace{-.4ex}\_\hspace{.7ex}$}S^1)^{\# 2}\#{\mathbb R}{\bf P}^3$
         & (14,84,140,70)   & 38 &    \\
         & [(15,87,144,72)] &    & 37 \\
 \addlinespace
 $(S^2\!\times\!S^1)^{\# 3}\#{\mathbb R}{\bf P}^3$,\, $(S^2\hbox{$\times\hspace{-1.62ex}\_\hspace{-.4ex}\_\hspace{.7ex}$}S^1)^{\# 3}\#{\mathbb R}{\bf P}^3$
         & (15,97,164,82) & 47 & 47 \\
 \addlinespace
 $(S^2\!\times\!S^1)^{\# 4}\#{\mathbb R}{\bf P}^3$,\, $(S^2\hbox{$\times\hspace{-1.62ex}\_\hspace{-.4ex}\_\hspace{.7ex}$}S^1)^{\# 4}\#{\mathbb R}{\bf P}^3$
         & (16,111,190,95) & 57 & 57 \\
 \addlinespace
 $(S^2\!\times\!S^1)^{\# 5}\#{\mathbb R}{\bf P}^3$,\, $(S^2\hbox{$\times\hspace{-1.62ex}\_\hspace{-.4ex}\_\hspace{.7ex}$}S^1)^{\# 5}\#{\mathbb R}{\bf P}^3$
         & (17,125,216,108) & 67 & 67 \\
 \addlinespace
 $(S^2\!\times\!S^1)\#L(3,1)$,\, $(S^2\hbox{$\times\hspace{-1.62ex}\_\hspace{-.4ex}\_\hspace{.7ex}$}S^1)\#L(3,1)$
         & (15,89,148,74)   & 46 &    \\  
         & [(16,92,152,76)] &    & 38 \\
 \addlinespace
 $(S^2\!\times\!S^1)^{\# 2}\#L(3,1)$,\, $(S^2\hbox{$\times\hspace{-1.62ex}\_\hspace{-.4ex}\_\hspace{.7ex}$}S^1)^{\# 2}\#L(3,1)$
         & (15,100,170,85)   & 50 &    \\
         & [(16,102,172,86)] &    & 48 \\
 \addlinespace
 $L(3,1)\#L(3,1)$
         & (16,113,194,97)   & 59 &    \\
         & [(17,115,196,98)] &    &    \\
         & [(20,126,212,106)]&    & 56 \\
 \addlinespace
 $L(3,1)\#-L(3,1)$
         & (16,115,198,99)   & 61 &    \\
         & [(17,115,196,98)] &    &    \\
         & [(20,126,212,106)]&    & 56 \\
 \addlinespace
 \addlinespace
\bottomrule
\end{tabular}
\end{table}

By our bistellar flip search it turned out that not always the triangulations
with fewest vertices have the smallest $g_2$ and therefore provide the best upper bound
on $\Gamma$.

\begin{thm}
The $3$-manifold ${\mathbb R}{\bf P}^{\,3}\#\,{\mathbb R}{\bf P}^{\,3}$ has (at least) two 
different minimal $g$-vectors. 
\end{thm}

\noindent
\textbf{Proof:}
The range of $f$-vectors for the projective space ${\mathbb R}{\bf P}^{\,3}$ is described by Walkup's Theorem~\ref{thm:Walkup} 
by $\Gamma^*=\Gamma =17$. The unique minimal triangulation of ${\mathbb R}{\bf P}^{\,3}$ has 
face vector $f=(11,51,80,40)$ 
and $g=(6,17)$. If we use this triangulation $K$ to form a triangulation
$K\# K$ of ${\mathbb R}{\bf P}^{\,3}\#\,{\mathbb R}{\bf P}^{\,3}$, then
$f(K\# K)=(18,96,156,78)$ and $g(K\# K)=(13,34)$. In particular,\,
$\Gamma({\mathbb R}{\bf P}^{\,3}\#\,{\mathbb R}{\bf P}^{\,3})\leq 34$.

On the other hand, by using bistellar flips, we obtained a triangulation of ${\mathbb R}{\bf P}^{\,3}\#\,{\mathbb R}{\bf P}^{\,3}$
with $f=(15,86,142,71)$ and $g=(10,36)$. This triangulation  showed up in our
enumeration of potential $g_2$-irreducible triangulations with up to $15$ vertices. Since there are
no $g_2$-irreducible $15$-vertex triangulations of ${\mathbb R}{\bf P}^{\,3}\#\,{\mathbb R}{\bf P}^{\,3}$ with 
$g_2<36$ and no $g_2$-irreducible triangulations of ${\mathbb R}{\bf P}^{\,3}\#\,{\mathbb R}{\bf P}^{\,3}$
with fewer vertices, the theorem follows.
\hfill$\Box$

\bigskip

Various other $3$-manifolds also seem to have non-unique minimal $g$-vectors. Candidates for such manifolds 
are listed in the Tables~\ref{tbl:small_lens_spaces}--\ref{tbl:connected_sums_manifolds} with additionally 
found $f$-vectors in brackets $[(f_0,f_1,f_2,f_3)]$. 
We use the refined simulated annealing process described earlier with a
restriction on $f_0$ during the cooling stage that it be at least $r$ greater than
the smallest known $f_0$ for the given $3$-manifold and a fixed $r=1$, $2$, $3$, or $4$.

For example, K\"uhnel and Lassmann \cite{KuehnelLassmann1984-3torus} described a neighborly $15$-vertex triangulation 
of the $3$-torus $T^3$ with $f=(15,105,180,90)$ and $g_2=55$. For this triangulation it is conjectured \cite{Lutz2005bpre}
that it is the unique vertex-minimal triangulation of $T^3$. By our search we found a triangulation of $T^3$ 
with $f=(16,108,184,92)$ and $g_2=54$. Therefore, $\Gamma(T^3)\leq 54$, which disproves $\Gamma(T^3)=55$
as was conjectured in \cite{Lutz2005bpre}. 

\bigskip

We finally have to describe how to obtain upper bounds on $\Gamma^*$.

\begin{lem}\label{lem:Hamiltonian}
Let $K$ be a neighborly triangulation of a $3$-manifold $M$ such that 
for some vertex $v$ of $K$ the link ${\rm Lk}\,v$ admits an Hamiltonian cycle 
going through all vertices of ${\rm Lk}\,v$. Then for every pair $(f_0,f_1)$ 
with $f_0\geq 0$ and
\begin{equation}
\binom{f_0}{2}\geq f_1\geq 4f_0 - 10 +g_2(K)
\end{equation}
there is a triangulation of $M$ with $f_0$ vertices and $f_1$ edges.
\end{lem}

\textbf{Proof}: The lemma is a reformulation of Walkup's Lemmas 7.3 from \cite{Walkup1970}.
The condition of the existence of an Hamiltonian cycle in the link of some vertex $v$
is equivalent to Walkup's condition of the existence of a ``spanning simple $3$-tree $T$'' 
whose $3$-simplices have the vertex $v$ in common.\hfill$\Box$

\bigskip

Let $K$ be a triangulation of some given $3$-manifold $M$. If $K$ is neighborly and fulfills
the requirement of Lemma~\ref{lem:Hamiltonian}, then $\Gamma^*(M)\leq g_2(K)$.
If $K$ is not neighborly we can try with bistellar flips to reach a neighborly triangulation $K'$ 
with the same number of vertices $f_0$ that fulfills the requirement of Lemma~\ref{lem:Hamiltonian}. 
Obviously, we then have $\Gamma^*(M)\leq g_2(K')$.

\begin{lem}\label{lem:upperbound}
Let $K$ be a triangulation with $f$-vector $f=(f_0,f_1,f_2,f_3)$ of a $3$-man\-i\-fold $M$ 
that can be connected to a neighborly triangulation $K'$ of $M$ with the same number of vertices $f_0$
via bistellar flips (with all intermediate triangulations also with $f_0$ vertices).
If $K'$ fulfills the requirement of Lemma~\ref{lem:Hamiltonian},
then
\begin{equation}
\Gamma^*(M)\,\leq\, \mbox{\rm max}\,\{\,g_2(K),\textstyle \binom{f_0-1}{2}-4(f_0-1)+10+1\,\}.
\end{equation}
\end{lem}

\textbf{Proof}: A neighborly triangulation $K''$ of $M$ with $f_0-1$ vertices,
if such a triangulation exists, would have $g_2(K'')=\binom{f_0-1}{2}-4(f_0-1)+10$.
If $g_2(K)\geq \binom{f_0-1}{2}-4(f_0-1)+10+1$, then there is no
triangulation $K'''$ of $M$ with fewer than $f_0$ vertices and $g_2(K''')\geq g_2(K)$. 
By the sequence of bistellar flips that connects $K$ with $K'$ we have
triangulations of $M$ with $f_0$ vertices for all integer values $g_2$ 
between $g_2(K)$ and $g_2(K')$. Applying subdivision operations $S$
to these triangulations guarantees the existence of triangulations of $M$
in the range $g_2(K)\leq g_2\leq g_2(K')$ for all number of vertices $n\geq f_0$.
In combination with $\Gamma^*(M)\leq g_2(K')$ it follows that $\Gamma^*(M)\leq g_2(K)$.

If $g_2(K)\leq\binom{f_0-1}{2}-4(f_0-1)+10$, then we can at least guarantee the upper bound
$\Gamma^*(M)\leq \binom{f_0-1}{2}-4(f_0-1)+10+1$.\hfill$\Box$

\bigskip
\bigskip

We used Lemmas~\ref{lem:Hamiltonian} and~\ref{lem:upperbound} to establish the $\Gamma^*$-bounds 
in the Tables~\ref{tbl:small_lens_spaces}--\ref{tbl:connected_sums_manifolds}.

\bigskip

\begin{thm}\label{thm:extensionWalkup}
A complete description of the set of all possible $f$-vectors is given by
\begin{itemize}
\addtolength{\itemsep}{-.5mm}
\item[{\rm (a)}] $\Gamma^*=\Gamma =0$\, for\, $S^3$,
\item[{\rm (b)}] $\Gamma^*=\Gamma =10$\, for\, $S^2\hbox{$\times\hspace{-1.62ex}\_\hspace{-.4ex}\_\hspace{.7ex}$}S^1$,

\item[{\rm (c)}] $\Gamma^*=11$\, and\, $\Gamma =10$\, for\,
  $S^2\!\times\!S^1$, where, with the exception $(9,36)$, all pairs $(f_0,f_1)$
  with $f_0\geq 0$ and $4f_0\leq f_1\leq\binom{f_0}{2}$ occur,

\item[{\rm (d)}] $\Gamma^*=\Gamma =17$\, for\, ${\mathbb R}{\bf P}^{\,3}$, 

\item[{\rm (e)}] $\Gamma^*=22$\, and\, $\Gamma =20$\, for\,
  $(S^2\!\times\!S^1)^{\# 2}$ and $(S^2\hbox{$\times\hspace{-1.62ex}\_\hspace{-.4ex}\_\hspace{.7ex}$}S^1)^{\# 2}$, 
  where, with the exceptions $(11,44)$ and $(11,45)$, all pairs $(f_0,f_1)$
  with $f_0\geq 0$ and $4f_0+10\leq f_1\leq\binom{f_0}{2}$ occur,

\item[{\rm (f)}] $\Gamma^*=\Gamma =10k$\, for\, $(S^2\!\times\!S^1)^{\# k}$ 
                 and $(S^2\hbox{$\times\hspace{-1.62ex}\_\hspace{-.4ex}\_\hspace{.7ex}$}S^1)^{\# k}$, $k=3,4,5,6,7,8,10,11,14$,
\addtolength{\itemsep}{.5mm}
\end{itemize}
while $\Gamma^*(M)\geq\Gamma (M)\geq 21$\, for all other $3$-manifolds $M$.
\end{thm}

\bigskip

\textbf{Proof:} Parts (a)--(d) follow from Theorem~\ref{thm:Walkup} of Walkup.
Since there are no $11$-vertex triangulations of $(S^2\!\times\!S^1)^{\# 2}$ 
and $(S^2\hbox{$\times\hspace{-1.62ex}\_\hspace{-.4ex}\_\hspace{.7ex}$}S^1)^{\# 2}$ \cite{SulankeLutz2006pre},
the lower bound (\ref{eq:NovikSwartz}) and the triangulations that were used to compute the respective 
upper bounds on $\Gamma$ and $\Gamma^*$ in Table~\ref{tbl:connected_sums_manifolds} together imply~(e).
The lower bound (\ref{eq:NovikSwartz}) and respective triangulations (cf.\ Table~\ref{tbl:connected_sums_manifolds})
imply~(f). 

Due to Theorem~\ref{thm:g2irred20} there is only one $g_2$-irreducible triangulation
of a $3$-manifold with $\Gamma\leq 20$, which is the unique $g_2$-irreducible triangulation
of ${\mathbb R}{\bf P}^{\,3}$. As observed in Section~\ref{sec:g2irreducible}, 
a triangulation $K$ which realizes the minimum $g_2$ for a particular $3$-manifold $M$ is either $g_2$-irreducible, 
or is of the form $K_1 \# K_2$ or $HK^\prime$, where the component triangulations realize their minimum $g_2$. 
In the case $K=K_1 \# K_2$ we have $g_2(K_1\# K_2) = g_2(K_1)+g_2(K_2)$ according to (\ref{eq:g2m1m2}),
whereas $g_2(HK^\prime)= g_2(K^\prime)+10$ in the case $K=HK^\prime$ according to (\ref{eq:g2hm3d}).
The $3$-manifolds obtained by adding a handle to one of the $3$-manifolds listed in (a)--(c) are listed
in (b), (c), and~(e). If $K^\prime$ triangulates ${\mathbb R}{\bf P}^{\,3}$ or a $3$-manifold $M$
with $\Gamma(M)\geq 18$, then $g_2(HK^\prime)\geq 27$. Similarly, all manifolds with triangulations
of the form $K=K_1 \# K_2$ (with one or both factors possibly triangulations of $S^3$)
are either listed in (a)--(e) or have $g_2(K_1\# K_2)\geq 27$. Thus $\Gamma^*(M)\geq\Gamma (M)\geq 21$\, 
for a $3$-man\-i\-fold~$M$ other than those from (a)--(e).\hfill$\Box$

\pagebreak

\begin{cor}
The following $3$-manifolds have unique componentwise minimal $f$-vec\-tors:

\vspace{-2mm}

\begin{itemize}
\addtolength{\itemsep}{-3.25mm}
\item $S^3$\, with\, $f=(5,10,10,5)$,
\item $S^2\hbox{$\times\hspace{-1.62ex}\_\hspace{-.4ex}\_\hspace{.7ex}$}S^1$\, with\, $f=(9,36,54,27)$,
\item $S^2\!\times\!S^1$\, with\, $f=(10,40,60,30)$,
\item ${\mathbb R}{\bf P}^3$\, with\, $f=(11,51,80,40)$,
\item $(S^2\!\times\!S^1)^{\# 2}$\, and\, $(S^2\hbox{$\times\hspace{-1.62ex}\_\hspace{-.4ex}\_\hspace{.7ex}$}S^1)^{\# 2}$\, with\, $f=(12,58,92,46)$,
\item $(S^2\!\times\!S^1)^{\# 3}$\, and\, $(S^2\hbox{$\times\hspace{-1.62ex}\_\hspace{-.4ex}\_\hspace{.7ex}$}S^1)^{\# 3}$\, with\, $f=(13,72,118,59)$,
\item $(S^2\!\times\!S^1)^{\# 4}$\, and\, $(S^2\hbox{$\times\hspace{-1.62ex}\_\hspace{-.4ex}\_\hspace{.7ex}$}S^1)^{\# 4}$\, with\, $f=(14,86,144,72)$,
\item $(S^2\!\times\!S^1)^{\# 5}$\, and\, $(S^2\hbox{$\times\hspace{-1.62ex}\_\hspace{-.4ex}\_\hspace{.7ex}$}S^1)^{\# 5}$\, with\, $f=(15,100,170,85)$,
\item $(S^2\!\times\!S^1)^{\# 6}$\, and\, $(S^2\hbox{$\times\hspace{-1.62ex}\_\hspace{-.4ex}\_\hspace{.7ex}$}S^1)^{\# 6}$\, with\, $f=(16,114,196,98)$,
\item $(S^2\!\times\!S^1)^{\# 7}$\, and\, $(S^2\hbox{$\times\hspace{-1.62ex}\_\hspace{-.4ex}\_\hspace{.7ex}$}S^1)^{\# 7}$\, with\, $f=(17,128,222,111)$,
\item $(S^2\!\times\!S^1)^{\# 8}$\, and\, $(S^2\hbox{$\times\hspace{-1.62ex}\_\hspace{-.4ex}\_\hspace{.7ex}$}S^1)^{\# 8}$\, with\, $f=(18,142,248,124)$,
\item $(S^2\!\times\!S^1)^{\# 10}$\, and\, $(S^2\hbox{$\times\hspace{-1.62ex}\_\hspace{-.4ex}\_\hspace{.7ex}$}S^1)^{\# 10}$\, with\, $f=(19,166,294,147)$,
\item $(S^2\!\times\!S^1)^{\# 11}$\, and\, $(S^2\hbox{$\times\hspace{-1.62ex}\_\hspace{-.4ex}\_\hspace{.7ex}$}S^1)^{\# 11}$\, with\, $f=(20,180,320,160)$,
\item $(S^2\!\times\!S^1)^{\# 14}$\, and\, $(S^2\hbox{$\times\hspace{-1.62ex}\_\hspace{-.4ex}\_\hspace{.7ex}$}S^1)^{\# 14}$\, with\, $f=(22,218,392,196)$.
\end{itemize}
\addtolength{\itemsep}{3.25mm}
\end{cor}

\begin{rem}
It is not known whether there is a $3$- or higher-dimensional manifold
that admits different minimal $f$-vectors.
\end{rem}

\section{\mathversion{bold}$\Gamma$-Values\mathversion{normal} and Matveev Complexity}
\label{sec:gamma_complexity}

In \cite{Matveev1990} Matveev introduced a notion of complexity for three-manifolds.  The main properties of $c(M),$ the {\em complexity} of $M,$ are

\begin{itemize}
  \item
   $c(M)$ is a nonnegative integer and
   for fixed $n,$ the number of irreducible closed three manifolds with $c(M) \le n$ is finite.
   \item
   $c(M_1 \# M_2) = c(M_1) + c(M_2).$
\end{itemize}
Since then, a significant amount of research has gone into understanding $c(M)$ and determining 
all closed $3$-manifolds with small complexity.  See the recent text \cite{Matveev2007} for details.  
In view of the fact that $\Gamma(M)$ also has the first two properties and Conjecture \ref{conj:connectedsum} 
is the third property, it is natural to consider the relationship between $c(M)$ and $\Gamma(M).$

For irreducible $3$-manifolds with $c(M)>0$ the complexity of $M$ is the number of tetrahedra 
needed in a pseudo-simplicial triangulation.  As noted earlier, two barycentric subdivisions 
produces a simplicial triangulation of $M$ with $24^2 c(M)$ facets.  As $g_2 \le h_2$ 
and the number of facets for any $K$ is $2+2h_1 + h_2,$ we see that $\Gamma(M) \le O(c(M)).$  
Conversely, if $M$ is irreducible, then any triangulation $K$ of $M$ which realizes $\Gamma(M)$ 
must be a $g_2$-irreducible triangulation. By Theorem \ref{thm:newwalkup}, the number of vertices 
in $K$ is bounded linearly by $\Gamma(M).$  Hence $h_1$ and the number of facets are linearly bounded 
by $\Gamma(M).$  So, $c(M) \le O(\Gamma(M)).$  It is easy to see that these estimates are very crude.

\begin{prob}
 What is the exact relationship between $c(M)$ and $\Gamma(M)$?
\end{prob}

\begin{table}
\small\centering
\defaultaddspace=0.2em
\caption{Irreducible $3$-manifolds sorted with respect to $\Gamma$-bounds up to $66$.}\label{tbl:sorted_by_gamma}
\begin{tabular}{@{}r@{\hspace{12mm}}l@{}}
\\\toprule
 \addlinespace
  Upper Bound for $\Gamma$ & Manifolds \\ \midrule
 \addlinespace
 \addlinespace
       0  &  $S^3$ \\
      17  &  ${\mathbb R}{\bf P}^3$ \\
      28  &  $L(3,1)$ \\
      38  &  $L(4,1)$, ${\mathbb R}{\bf P}^{\,2}\!\times S^1$ \\
      40  &  $L(5,2)$, $P_2=S^3/Q$ \\
      47  &  $L(5,1)$, $P_3$ \\
      50  &  $L(7,2)$, $S^3/T^*$\\
      51  &  $L(8,3)$ \\
      52  &  $S^3/I^*=\Sigma(2,3,5)$ \\
      54  &  $P_4$, $T^3$ \\
      55  &  $S^3/O^*$ \\
      56  &  $L(6,1)$, $B_2$ \\
      57  &  $G_5$ \\
      58  &  $L(9,2)$ \\
      59  &  $L(10,3)$, $G_3$, $B_4$, $\{Oo,1\mid 1\}$ \\
      60  &  $K\times S^1$, $G_2$ \\
      61  &  $G_4$, $B_3$ \\
      62  &  $\{Oo,1\mid 2\}$ \\
      63  &  $\Sigma(2,3,7)$ \\
      64  &  $P_5$ \\
      65  &  $L(7,1)$ \\
      66  &  $G_6$ \\
 \addlinespace
 \addlinespace
\bottomrule
\end{tabular}
\end{table}

\begin{table}
\small\centering
\defaultaddspace=0.2em
\caption{Irreducible orientable $3$-manifolds of complexity up to $3$.}\label{tbl:sorted_by_complexity}
\begin{tabular}{@{}r@{\hspace{12mm}}l@{}}
\\\toprule
 \addlinespace
   Complexity & Manifolds \\ \midrule
 \addlinespace
 \addlinespace
       0  &  $S^3$, ${\mathbb R}{\bf P}^3$, $L(3,1)$ \\
 \addlinespace
       1  &  $L(4,1)$, $L(5,2)$ \\
 \addlinespace
       2  &  $L(5,1)$, $L(7,2)$, $L(8,3)$, $P_2=S^3/Q$ \\
 \addlinespace
       3  &  $L(6,1)$, $L(9,2)$, $L(10,3)$, $L(11,3)$, $L(12,5)$, $L(13,5)$, $P_3$ \\
 \addlinespace
 \addlinespace
\bottomrule
\end{tabular}
\end{table}

Tables \ref{tbl:sorted_by_gamma} and \ref{tbl:sorted_by_complexity} (which is taken from \cite{Matveev1998pre}) 
show that, at least for manifolds of small complexity for which we have some data, the rank ordering of the 
two measures of complexity are very close. Indeed, we believe that the three irreducible $3$-manifolds 
with smallest $\Gamma$ are $\Gamma(S^3) = 0$, $\Gamma(\mathbb{R}P^3) = 17$ and $\Gamma(L(3,1)) = 28.$

\subsection*{Acknowledgment}

Some of this work was done while the  first and third authors were at  the special 
semester on ``Combinatorics of Polytopes and Complexes: Relations with Topology and Algebra" 
(spring and summer 2007) at the Institute for Advanced Studies in 
Jerusalem. We are grateful to IAS for their hospitality and especially to Gil Kalai 
for organizing this semester. We also thank the anonymous referees for helpful comments.

\bibliography{}

\end{document}